\documentclass{article}

\usepackage[includeheadfoot,
            bindingoffset = 0mm,
            inner = 25mm,
            top = 20mm,
            outer = 25mm,
            bottom = 10mm,
            paperwidth = 210mm,
            paperheight = 297mm,
            ]{geometry}
\usepackage[margin={2.0cm,0cm},oneside,labelfont={sf,bf},singlelinecheck=false]{caption}

\usepackage{graphicx}

\usepackage[fleqn]{amsmath}
\usepackage{accents}
\usepackage{amssymb}
\usepackage{stmaryrd}
\SetSymbolFont{stmry}{bold}{U}{stmry}{m}{n}

\usepackage{enumitem}
\usepackage{cite}
\usepackage[perpage,multiple]{footmisc}
\usepackage{hyphenat}
\usepackage[english]{babel}
\usepackage[section]{placeins}
\usepackage{upgreek}
\usepackage{color}
\usepackage{multirow}

\usepackage{amsthm}
\theoremstyle{definition}

\usepackage{titlesec}
\titleformat{\section}[hang]{\Large\bfseries\raggedright\sffamily}{\thesection}{1em}{}
\titleformat{\subsection}[hang]{\large\bfseries\raggedright\sffamily}{\thesubsection}{1em}{}
\titleformat{\subsubsection}[hang]{\normalsize\bfseries\raggedright\sffamily}{\thesubsubsection}{1em}{}
\usepackage{abstract}

\usepackage{authblk}

\newcommand{\vn}[1]{\ensuremath{ \boldsymbol{#1} }}

\newcommand{\transp}{\ensuremath{ ^\mathrm{T} }}

\newcommand{\dif}{\ensuremath{ \mathrm{d} }}

\begin{document}

\title{ \huge\bfseries\sffamily Construction of the Kolmogorov-Arnold representation using the Newton-Kaczmarz method }

\author[1,3]{Michael Poluektov}
\author[2]{Andrew Polar}
\affil[1]{Department of Mathematical Sciences and Computational Physics, School of Science and Engineering, University of Dundee, Dundee DD1 4HN, UK}
\affil[2]{Independent software consultant, Duluth, GA, USA}
\affil[3]{Corresponding author, email: mpoluektov001@dundee.ac.uk}

\date{ \huge\normalfont\sffamily DRAFT: \today }

\maketitle

\setlength{\absleftindent}{2.0cm}
\setlength{\absrightindent}{2.0cm}
\setlength{\absparindent}{0em}
\begin{abstract}
It is known that any continuous multivariate function can be represented exactly by a composition functions of a single variable --- the so-called Kolmogorov-Arnold representation. It can be a convenient tool for tasks where it is required to obtain a predictive model that maps some vector input of a black box system into a scalar output. In this case, the representation may not be exact, and it is more correct to refer to such structure as the Kolmogorov-Arnold model (or, as more recently popularised, `network'). Construction of such model based on the recorded input-output data is a challenging task. In the present paper, it is suggested to decompose the underlying functions of the representation into continuous basis functions and parameters. It is then proposed to find the parameters using the Newton-Kaczmarz method for solving systems of non-linear equations. The algorithm is then modified to support parallelisation. The paper demonstrates that such approach is also an excellent tool for data-driven solution of partial differential equations. Numerical examples show that for the considered model, the Newton-Kaczmarz method for parameter estimation is efficient and more robust with respect to the section of the initial guess than the straightforward application of the Gauss-Newton method. Finally, the Kolmogorov-Arnold model is compared to the MATLAB's built-in neural networks on a relatively large-scale problem ($25$ inputs, datasets of $10$ million records), significantly outperforming the multilayer perceptrons (MLPs) in this particular problem ($4$--$10$ minutes vs. $4$--$8$ hours of training time, as well as higher accuracy, lower CPU usage, and smaller memory footprint). 
Kolmogorov-Arnold representation, Kolmogorov-Arnold networks, ridge function, \\
Kaczmarz method, Newton-Kaczmarz method.
\end{abstract}

\section{Introduction}
\label{sec:intro}

One of the typical data science tasks is constructing a predictive model that maps some vector input into a scalar output. For example, this can be a regression model that establishes a relationship between the atomistic-scale structure of matter and its properties \cite{Musil2021}. Such model can be viewed as a multivariate function. In 1950s, Andrei Kolmogorov and Vladimir Arnold showed that any continuous multivariate function can be represented exactly by a composition of functions of a single variable \cite{Arnold1957,Kolmogorov1956}. Thus, the Kolmogorov-Arnold representation is convenient tool for data modelling tasks, replacing the problem of construction of a multivariate function in its entirety by the task of building of a set of functions of a single variable. In 1980s, an interpretation of the Kolmogorov-Arnold representation as a neural network was suggested \cite{HechtNielsen1987}, creating two directions --- research on the exact representation and on its approximate forms, which can be called \emph{models} or \emph{networks} to distinguish from the exact form.

As for the exact representation, the original Kolmogorov's proof was not constructive and did not provide insights into choosing the underlying functions of the representation. Subsequent refinements and generalisations of the theorem included unification of the outer functions \cite{Lorentz1962}, reduction of the number of the inner functions by introducing shifts \cite{Sprecher1965,Sprecher1972}, and improvement of the smoothness of the inner functions \cite{Fridman67}. More recent works focused on development of a constructive proof of the theorem with convergent algorithms for determining the underlying functions \cite{Sprecher1996,Sprecher1997,Koeppen2002,Braun2009,Actor2018}, and on generalisation of the theorem to discontinuous bounded and unbounded multivariate functions \cite{Ismayilova2024}.

The main effort to construct approximate forms of the representation started in 1990s, e.g. \cite{Nees1994}. Following ideas from neural networks, it has been proposed to approximate the underlying functions using sigmoids \cite{Kurkova1992}, which further bridges the representation and neural networks. The latter connection has been explored across a range of papers, including recent \cite{Montanelli2020,SchmidtHieber2020}, which focus on construction of deep ReLU networks that imitate the structure of the Kolmogorov-Arnold representation. From 2000s, the model has been used for machine-learning applications, such as image processing \cite{Bryant2008,Liu2015}.

A remarkable achievement of early 2000s should be especially highlighted --- it was proposed to take the underlying functions of the Kolmogorov-Arnold model in a form of cubic splines \cite{Igelnik2003}. An independent similar spline-based model architecture has been proposed a year later \cite{Coppejans2004}. Another spline-based model has been suggested in a recent preprint \cite{Deventer2022}. The applications of the spline-based Kolmogorov-Arnold model included image decomposition \cite{Leni10}.

For applications where smoothness of the underlying functions of the Kolmogorov-Arnold model is not essential, piecewise-linear approximation can be a good choice. Recently, a lightweight iterative algorithm for construction of the underlying functions limited to such class has been proposed by the authors of the present paper \cite{Polar2021}; however, although well-tested computationally, the algorithm was given without the convergence proof.

The present paper proposes to decompose the underlying functions of the Kolmogorov-Arnold model into continuous basis functions and parameters. It is then shown that having the input-output data, the parameters can be found efficiently using the Newton-Kaczmarz method for solving systems of non-linear equations, irrespective of the choice of the basis functions. The method is locally convergent. Finally, it is shown that the Kolmogorov-Arnold model is an efficient tool for data-driven solution of non-linear partial differential equations (PDEs).

It should be noted that although the Kolmogorov-Arnold representation/model, has been extensively researched since 1950s, it gained significant media attention in May 2024 with the publication of a preprint by Z. Liu \emph{et al}. \cite{Liu2024}, where the approximate form was named ``Kolmogorov-Arnold Network''. Although news items, blogs, and social media posts presented it as a breakthrough, this work was largely a retelling of the preceding ideas --- it suggested to approximate the underlying functions using splines (and possibly other functions) and to introduce more than two layers --- both ideas were not new. As mentioned above, the former was first proposed in 2003 by Igelnik and Parikh, whose paper was even entitled ``Kolmogorov's Spline Network'', and whose developments were used for practical applications in late 2000s by Leni \cite{Leni10}. Efficient use of arbitrary basis functions that are decoupled from the model's training algorithm was shown in the first preprint of the present paper\footnote{https://arxiv.org/abs/2305.08194v1}, posted online in May 2023, one year prior to \cite{Liu2024}. The idea of introducing more than two layers (deep network) was at least formulated by the authors of the present paper in \cite{Polar2021} and then computationally tested together with an uncertainty quantification algorithm in a preprint\footnote{https://arxiv.org/abs/2104.01714v3}\footnote{The interested readers should look specifically at version $3$, since the scope of the subsequent versions of the paper has been shifted, excluding multi-layer networks.} in 2021; moreover, it is highly likely that it was first formulated much earlier by pioneering works on the Kolmogorov-Arnold model from 1990s and early 2000s but did not gain traction because of the lack of visible advantages (compared to the `classical' two-layer version) on small-scale problems. The novel aspect of \cite{Liu2024} was a model reduction technique consisting in switching off least involved functions of the composition. After May 2024, there has been an extremely rapid growth of preprints on the Kolmogorov-Arnold model (now called KAN), but an extensive overview of the latter is out of the scope of the present paper.

Finally, it should be emphasised that the authors do not intend to convince the reader that the Kolmogorov-Arnold model/network is the best possible machine learning model. Question whether this model can compete with the state-of-the-art machine-learning models on real large-scale datasets will remain open until a large range of comprehensive studies by a vast number of groups using a variety of implementations is performed. Some steps are already made --- for example, some recent studies already apply the authors' open-source implementation, developed for the present paper, to industrial problems \cite{Sulaiman2024,Sulaiman2024a}. The Kolmogorov-Arnold model (i.e. the approximate form of the representation) exists in literature for at least 30 years in various forms and has been used in practice, which motivates developing fast and reliable training methods for it --- the scope of the present paper. The comparisons to neural networks presented below should not be interpreted as an ultimate evidence in any favour.

The Kolmogorov-Arnold model is related to some simpler models, in particular, to the Urysohn/GAM model. The parameter estimation techniques for these models are related. Therefore, the main part of the present paper starts with the description of the Urysohn/GAM model and its construction using the data.

\section{Models}
\label{sec:models}

Since the Kolmogorov-Arnold representation/model has been researched by the applied maths community, for historic reasons, the present paper will use the applied maths terminology: \emph{inputs} are \emph{features}, \emph{outputs} are \emph{targets}, \emph{data records} are \emph{samples}, and \emph{passes} are \emph{epochs}. 

\subsection{The Urysohn/GAM model}
\label{sec:modelsUry}

Given set of data records $\left\lbrace \vn{x}_i, y_i \right\rbrace$, $i \in \left\lbrace 1, \ldots, N \right\rbrace$, where $\vn{x}_i \in \mathbb{R}^m$ is the input of the $i$-th record, $y_i \in \mathbb{R}$ is the output of the $i$-th record, and $N$ is the number of records, the following model can be introduced:
\begin{equation}
  \tilde{y}_i = \sum_{j=1}^{m} g_j \left(\vn{x}_{i,j}\right) ,
  \label{eq:UrysDiscr}
\end{equation} 
where $x_{i,j} \in \left[x_\mathrm{min}, x_\mathrm{max}\right]$ denotes the $j$-th component of vector $\vn{x}_i$, and $\tilde{y}_i$ is the calculated model output of the $i$-th record. Functions $g_j: \left[x_\mathrm{min}, x_\mathrm{max}\right] \to \mathbb{R}$ are continuous almost everywhere. It is important to emphasise the distinction between actual (recorded) output $y_i$ and model output $\tilde{y}_i$. 

Model \eqref{eq:UrysDiscr} can have different names. It was originally written for time-series data and most research on it has been done within the fields of automatic control and signal processing, as summarised in \cite{Poluektov2020} and references therein. In this case, the model maps one time series into another, and it is called the discrete Urysohn operator.

Model \eqref{eq:UrysDiscr} is also referred to as the generalised additive model (GAM) in literature \cite{Hastie1990}. The latter has been introduced and researched independently to the time-series case.

\subsection{The Kolmogorov-Arnold model}
\label{sec:KA}

In \cite{Arnold1957,Kolmogorov1956}, it has been shown that any continuous multivariate function can be represented by a composition of functions of a single variable. More precisely, any continuous function $F: \left[0,1\right]^m \to \mathbb{R}$ can be exactly represented in the following way:
\begin{equation}
  F\left( x_1, x_2, \ldots, x_m \right) = \sum_{k=1}^{2m+1} \varPhi_k \left( \sum_{j=1}^{m} f_{kj} \left(x_j\right) \right) ,
  \label{eq:Kolmogorov}
\end{equation}
where $f_{kj}: \left[0,1\right] \to \mathbb{R}$ and $\varPhi_k: \mathbb{R} \to \mathbb{R}$ are continuous functions. This decomposition can be used for the tasks of data modelling, where the output is a continuous function of the inputs. Given set of data records $\left\lbrace \vn{x}_i, y_i \right\rbrace$, intermediate variable $\vn{\theta}_i \in \mathbb{R}^{d}$ is introduced for each data record, and the Kolmogorov-Arnold model for the input-output relationship is then written as
\begin{equation}
  \tilde{y}_i = \sum_{k=1}^{d} \varPhi_k \left( \theta_{i,k} \right) , \quad
  \theta_{i,k} = \sum_{j=1}^{m} f_{kj} \left( x_{i,j} \right) , \quad
  k \in \left\lbrace 1, \ldots, d \right\rbrace ,
  \label{eq:KAU}
\end{equation} 
where $\theta_{i,k}$ denotes the $k$-th component of vector $\vn{\theta}_i$ and $d$ is the number of addends in the model; $d=2m+1$ for the `classical' representation, aligned with equation \eqref{eq:Kolmogorov}. For consistency of the notation, the domains of functions $f_{kj}$ have been extended to $\left[x_\mathrm{min}, x_\mathrm{max}\right]$ in equation \eqref{eq:KAU}, which can always be done by rescaling of the inputs. It should be mentioned that the Kolmogorov-Arnold model can also be constructed of discontinuous functions \cite{Ismayilova2024}; however, no practical parameter estimation algorithms for a model with discontinuous or unbounded functions have been proposed.

This model can also be interpreted as a particular tree (or a particular `network') of the Urysohn/GAM models, containing $1$ `root' model with $d$ functions and $d$ `branch' models with $m$ functions each. Intermediate variable $\vn{\theta}_i$ can then be interpreted as a hidden variable between two layers.

It is useful to note that the Kolmogorov-Arnold representation can also be viewed as a generalisation of the so-called ridge function \cite{Ismailov2020}:
\begin{equation}
  F\left( x_1, x_2, \ldots, x_m \right) = \varPhi\left( \sum_{j=1}^{m} c_j  x_j \right) ,
  \label{eq:ridge}
\end{equation}
where $c_j \in \mathbb{R}$ are some constants and $\varPhi: \mathbb{R} \to \mathbb{R}$ is some function. The ridge function models are used in some areas of data modelling and are even suggested as an approximate replacement of the Kolmogorov-Arnold model and neural networks \cite{Ismailov2020}.

\section{Discretisation and parameter estimation}
\label{sec:discr}

\subsection{Building the Urysohn/GAM model}
\label{sec:Uident}

To handle the models computationally, the underlying functions must be represented in some discretised form. For the Urysohn/GAM model, one option is to use the following decomposition of functions $g_j$:
\begin{equation}
  g_j\left(x\right) = \sum_{p=1}^{n} U_{jp} \phi_p\left(x\right) ,
  \label{eq:PWlinear}
\end{equation}
where $U_{jp}$ are the model parameters and $\phi_p: \left[x_\mathrm{min}, x_\mathrm{max}\right] \to \mathbb{R}$ are the basis functions with compact support. From equations \eqref{eq:UrysDiscr} and \eqref{eq:PWlinear}, it can be seen that parameters $U_{jp}$ can also be understood as elements of $m$-by-$n$ matrix $U$.

Training of the Urysohn/GAM model consists in finding parameters $U_{jp}$ given set of data records $\left\lbrace \vn{x}_i, y_i \right\rbrace$, such that the discrepancy between $\tilde{y}_i$ and $y_i$ is minimised. This problem and its properties (e.g. the non-uniqueness of the solution) have been studied in detail in \cite{Poluektov2020} for piecewise-linear and piecewise-constant basis functions. In the latter, the time-series form of the model was used; therefore, it is useful to rewrite the main result of \cite{Poluektov2020} for arbitrary basis functions and for the case of a set of data records. 

It can be seen that $\tilde{y}_i$ is a linear function of the parameters (but not a linear function of the inputs). Therefore, it is possible to assemble a system of linear equations with respect to $U_{jp}$ and use it to reformulate the parameter estimation problem. Thus, the model parameters are rearranged into vector $\vn{Z}$, the values of the basis functions are rearranged into matrix $\vn{M}$, rows of which are denoted as $\vn{M}_i$, and the outputs are rearranged into vector $\vn{Y}$:
\begin{align*}
  &\vn{Z} = \begin{bmatrix}
    U_{11} & \!\!\ldots\!\! & U_{1n} & U_{21} & \!\!\ldots\!\! & U_{2n} & \ldots & U_{m1} & \!\!\ldots\!\! & U_{mn}
  \end{bmatrix}\transp , \\
  &\vn{M}_i = \begin{bmatrix}
    \phi_1\left(x_{i,1}\right) & \!\!\ldots\!\! & \phi_n\left(x_{i,1}\right) & \phi_1\left(x_{i,2}\right) & \!\!\ldots\!\! & \phi_n\left(x_{i,2}\right) & \ldots & \phi_1\left(x_{i,m}\right) & \!\!\ldots\!\! & \phi_n\left(x_{i,m}\right)
  \end{bmatrix} , \\
  &\vn{Y} = \begin{bmatrix} y_1 & \!\!\ldots\!\! & y_N \end{bmatrix}\transp .
\end{align*}
This leads to the following minimisation problem:
\begin{equation}
  \vn{Z}^* = \underset{\vn{Z}}{\operatorname{argmin}} \left\|\vn{Y} - \vn{M} \vn{Z}\right\|^2 ,
  \label{eq:UryMin}
\end{equation}
where $\vn{Z}^*$ is the sought solution and $\left\|\cdot\right\|$ denotes the standard $l^2$-norm of a vector. In \cite{Poluektov2020}, it has been suggested to find an approximate solution of this problem iteratively row-by-row (or record-by-record) using the projection descend method \cite{Kaczmarz1937,Tewarson1969}, also known as the Kaczmarz method, which represents the following sequence:
\begin{equation}
  \vn{Z}^{q+1} = \vn{Z}^q + \mu\frac{y_i - \vn{M}_i \vn{Z}^q}{\left\| \vn{M}_i \right\|^2} \vn{M}_i^\mathrm{T} ,
  \label{eq:KaczGen}
\end{equation}
where $\mu \in \left(0,1\right]$ is the numerical damping parameter and $\vn{Z}^q$ is the approximation of the solution at iteration $q$. Index $i$ changes each iteration. The method requires initial guess $\vn{Z}^0$ and continues until some convergence criteria are reached. 

A number of important remarks should be made briefly regarding the method above (the details are given in \cite{Poluektov2020} and references therein). First, the method considers system of linear equations $\vn{M} \vn{Z} = \vn{Y}$. For an underdetermined system, the method converges to a solution that is closest to the initial guess in terms of the $l^2$-norm. For an overdetermined system, it produces a solution that oscillates in a region around the true solution of minimisation problem \eqref{eq:UryMin}. The size of this region decreases with decreasing $\mu$. Second, although an approximate solution is obtained, there is a significant advantage of the method --- assembling of matrix $\vn{M}$ and keeping it in the memory is not required, as the method iterates row-by-row. Furthermore, in the case of the Urysohn/GAM model, since functions $\phi_p$ have compact support, matrix $\vn{M}$ is sparse\footnote{For example, for piecewise-linear basis functions given by equation \eqref{eq:nodal}, it is easy to see that for any $x_{i,j}$, there will be at most two non-zero $\phi_p(x_{i,j})$ out of $n$. Hence, $\vn{M}^i$ will contain at most $2m$ non-zero elements, and the fraction of the non-zero elements of matrix $\vn{M}$ is at most $2/n$.}; therefore, only a small subset of elements of $\vn{Z}^q$ is updated in equation \eqref{eq:KaczGen}, which is computationally cheap. Third, a consequence of the sparsity of $\vn{M}$ is that the projection descend method converges relatively fast, as the rows of $\vn{M}$ are either orthogonal or close to orthogonal (the detailed analysis of the convergence of the projection descend method can be found in e.g. \cite{Popa2018} and references therein).

\subsection{Building the Kolmogorov-Arnold model}
\label{sec:KAident}

As in the previous subsection, the underlying functions of the model are decomposed as
\begin{align}
  &f_{kj}\left(x\right) = \sum_{p=1}^{n} H_{kjp} \phi_p\left(x\right) , 
  \label{eq:KAfunc1} \\
  &\varPhi_k\left(t\right) = \sum_{l=1}^{s} G_{kl} \psi_l\left(t\right) ,
  \label{eq:KAfunc2}
\end{align}
where $H_{kjp}$ and $G_{kl}$ are the model parameters, $\phi_p: \left[x_\mathrm{min}, x_\mathrm{max}\right] \to \mathbb{R}$ and $\psi_l: \left[t_\mathrm{min}, t_\mathrm{max}\right] \to \mathbb{R}$ are the basis functions, and $t_\mathrm{min}, t_\mathrm{max} \in \mathbb{R}$. Functions $\phi_p$ and $\psi_l$ must be continuous, and functions $\psi_l$ must be differentiable everywhere except at a finite number of points. Integers $n$ and $s$ are the numbers of the basis functions used for the inner and the outer functions of the model, respectively.

Training of the Kolmogorov-Arnold model consists in finding parameters $H_{kjp}$ and $G_{kl}$ given set of data records $\left\lbrace \vn{x}_i, y_i \right\rbrace$, such that the discrepancy between $\tilde{y}_i$ and $y_i$ is minimised. As above, this problem can be written as minimisation problem
\begin{equation}
  \left\lbrace H_{kjp}^* , G_{kl}^* \right\rbrace = \underset{\vn{H},\vn{G}}{\operatorname{argmin}} \sum_{i=1}^{N} \left(\tilde{y}_i - y_i\right)^2 ,
  \label{eq:KAMin}
\end{equation}
where $H_{kjp}^*$ and $G_{kl}^*$ are the sought parameters. In contrast to the Urysohn/GAM model, output $\tilde{y}_i$ is now a non-linear function of the model parameters. 

An approximate solution of problem \eqref{eq:KAMin} can be found using the Newton-Kaczmarz method \cite{Meyn1983,MARTINEZ1986,MARTINEZ1986a} for solving systems of non-linear equations. As outlined in appendix \ref{sec:NKmethod}, this method can be understood as an iterative sequence of alternating between equations' linearisation and application of one step of the Kaczmarz method. 

The system of non-linear equations with respect to the model parameters is written as
\begin{equation}
  \tilde{y}_i - y_i = 0 , \quad\quad
  i \in \left\lbrace 1, \ldots, N \right\rbrace .
  \label{eq:syst}
\end{equation}
The solution is found iteratively. The approximations of the model parameters at iteration $q$ are denoted as $H_{kjp}^q$ and $G_{kl}^q$. It is assumed that initial guesses $H_{kjp}^0$ and $G_{kl}^0$ for the parameters are made. Writing the Newton-Kaczmarz method (summarised in appendix \ref{sec:NKmethod}) for system \eqref{eq:syst} leads to 
\begin{align}
  &H_{kjp}^{q+1} = H_{kjp}^q - \mu \left(\tilde{y}_i - y_i\right)\zeta^{-1} \frac{\partial \tilde{y}}{\partial H_{kjp}} , \label{eq:iter0} \\
  &G_{kl}^{q+1} = G_{kl}^q - \mu \left(\tilde{y}_i - y_i\right)\zeta^{-1} \frac{\partial \tilde{y}}{\partial G_{kl}} , \label{eq:iter} 
\end{align}
where $\mu \in \left(0,1\right]$ is the numerical damping parameter. Variable $\zeta$ is the squared norm of the gradient of $\tilde{y}$ with respect to the parameters:
\begin{equation}
  \zeta = \sum_{k=1}^{d} \sum_{j=1}^{m} \sum_{p=1}^{n} \left(\frac{\partial \tilde{y}}{\partial H_{kjp}}\right)^2 + \sum_{k=1}^{d} \sum_{l=1}^{s} \left(\frac{\partial \tilde{y}}{\partial G_{kl}}\right)^2 .
\end{equation}
The derivatives of $\tilde{y}$ are summarised in appendix \ref{sec:derv}. It is implied that the values of the derivatives are calculated at $\vn{x}_i$ and using the parameters at iteration $q$. As previously, index $i$ changes each iteration. The iterations are performed until some convergence criteria are reached, e.g. $\left|\tilde{y}_i - y_i\right|$ is sufficiently small for sufficiently large number of iterations (the exact thresholds are decided by the user). 

As for the case of the Urysohn/GAM model, system \eqref{eq:syst} can be underdetermined or overdetermined. An example of the latter is a system with noise-affected recorded output $y_i$. As in section \ref{sec:Uident}, the method then gives an approximate solution, which oscillates in a region around the true solution of minimisation problem \eqref{eq:KAMin} with the size of this region depending on $\mu$.

It is important to note that the original Newton-Kaczmarz method is formulated for systems of non-linear equations where functions are continuously differentiable, as outlined in appendix \ref{sec:NKmethod}. Thus, if functions $\psi_l$ are continuously differentiable, then the algorithm is locally convergent, as it is the direct application of the Newton-Kaczmarz method to system \eqref{eq:syst}, which makes it distinct from the algorithm proposed previously in \cite{Polar2021} that was given without the convergence proof.

If the derivatives of functions $\psi_l$ are discontinuous only at points belonging to a finite set, iterations at which $\theta_{i,k}$ falls exactly at a point of discontinuity can be simply skipped. In practice, this does not impact the convergence of the method. The relaxation of the continuity requirements allows using piecewise-linear basis functions, which is convenient from a practical point of view.

\subsubsection{Piecewise-linear basis functions}

In practice, it is convenient to use the piecewise-linear basis functions because of simplicity. To define the functions, equally-spaced nodes are introduced as
\begin{equation}
  x_p = x_\mathrm{min} + \left(p-1\right) \Delta x , 
    \quad\quad \Delta x = \frac{x_\mathrm{max} - x_\mathrm{min}}{n-1} , 
    \quad\quad p \in \left\lbrace 1, \ldots, n \right\rbrace .
  \label{eq:nodes}
\end{equation}
The same is done for the intermediate variable:
\begin{equation}
  t_l = t_\mathrm{min} + \left(l-1\right) \Delta t , 
    \quad\quad \Delta t = \frac{t_\mathrm{max} - t_\mathrm{min}}{s-1} ,
    \quad\quad l \in \left\lbrace 1, \ldots, s \right\rbrace .
  \label{eq:nodesInter}
\end{equation}
As introduced previously, $n$ and $s$ are the numbers of the inner and the outer basis functions, respectively. It should be noted that in many problems, the optimal choice of nodes is not equally-spaced, and redistributing the nodes adaptively can improve the accuracy of the model, however, leading to a more complex code/implementation. For convenience, the nodes are combined in sets $P = \left\lbrace x_p \right\rbrace$ and $Q = \left\lbrace t_l \right\rbrace$. The basis functions are continuous, are linear between the nodes, are exactly $1$ at a given node and are exactly $0$ at all other nodes:
\begin{equation}
  \phi_p\left(x\right) = \begin{cases}
    0 , & \text{if } x \in P \text{ and } x \neq x_p , \\
    1 , & \text{if } x = x_p , \\
    \text{linear} , & \text{otherwise} ,
    \end{cases} \quad\quad\quad
  \psi_l\left(t\right) = \begin{cases}
    0 , & \text{if } t \in Q \text{ and } t \neq t_l , \\
    1 , & \text{if } t = t_l , \\
    \text{linear} , & \text{otherwise} .
    \end{cases}
  \label{eq:nodal}
\end{equation}

\subsubsection{Gaussian basis functions}

A simple choice of smooth basis functions is a set of Gaussian functions. Using equally-spaced nodes introduced above, these functions can be formally written as 
\begin{align}
  &\phi_p\left(x\right) = \exp\left(-\frac{\gamma}{{\Delta x}^2}\left(x-x_p\right)^2\right) , \\
  &\psi_l\left(t\right) = \exp\left(-\frac{\gamma}{{\Delta t}^2}\left(t-t_l\right)^2\right) , \label{eq:gaussbf}
\end{align}
where $\gamma$ is the decay-rate parameter.

\subsubsection{Spline basis functions}

Another convenient choice of basis functions is a set of spline functions. This choice is common in numerical approximation and forms the foundation of techniques such as isogeometric analysis \cite{Cottrell2009}. Splines have been previously used to approximate functions of the Kolmogorov-Arnold model \cite{Deventer2022} and of multi-layer perceptrons \cite{Lane1990}.

The classical choice for 1D interpolation problems is cubic splines. Thus, the basis functions can be formally written as twice continuously differentiable functions such that
\begin{equation}
  \phi_p\left(x\right) = \begin{cases}
    0 , & \text{if } x \in P \text{ and } x \neq x_p , \\
    1 , & \text{if } x = x_p , \\
    \text{cubic} , & \text{otherwise} ,
    \end{cases} \quad\quad\quad
  \psi_l\left(t\right) = \begin{cases}
    0 , & \text{if } t \in Q \text{ and } t \neq t_l , \\
    1 , & \text{if } t = t_l , \\
    \text{cubic} , & \text{otherwise} .
    \end{cases}
  \label{eq:splinebf}
\end{equation}
The latter must be supplemented with the end conditions, such as \emph{not-a-knot} or \emph{natural}.

\subsubsection{Initial guess for the parameters and the model complexity}

In general, initial approximations $H_{kjp}^0$ and $G_{kl}^0$ can be arbitrary. However, the outputs of the inner functions constitute the intermediate variable, the components of which are discretised according to equation \eqref{eq:nodesInter}. With `bad' initial approximations of the parameters, the intermediate variable might significantly drift during the parameter estimation process. In the previous method \cite{Polar2021}, it was suggested to update limits $t_\mathrm{min}$ and $t_\mathrm{max}$ of the intermediate variable and its nodal positions $t_l$. A better approach is to initialise the parameters in such a way that this updating procedure is not necessary.

One such way is to take $t_\mathrm{min} = y_\mathrm{min}$ and $t_\mathrm{max} = y_\mathrm{max}$, where $y_\mathrm{min}$ and $y_\mathrm{max}$ are the minimum and the maximum values of output $y_i$, respectively. The model parameters are then initialised as uniformly-distributed random numbers:
\begin{equation}
  H_{kjp}^0 \sim \operatorname{unif}\left(\frac{y_\mathrm{min}}{m},\frac{y_\mathrm{max}}{m}\right) , \quad\quad 
  G_{kl}^0 \sim \operatorname{unif}\left(\frac{y_\mathrm{min}}{d},\frac{y_\mathrm{max}}{d}\right) .
  \label{eq:initHG}
\end{equation}

The complexity of the Kolmogorov-Arnold model is determined by the numbers of the basis functions, $n$ and $s$, and by the number of addends, $d$. The total number of the model parameters is $d\left(nm+s\right)$. The `full' (or `classical') model corresponds to $d=2m+1$ that aligns with the Kolmogorov-Arnold representation theorem. Depending on the complexity of the modelled system and the desired accuracy, hyper-parameters $n,s,d$ should be fine-tuned using the training and the selection (parameter optimisation) datasets, with subsequent assessment using the validation dataset.

\subsubsection{Parallelisation/vectorisation of the algorithm}
\label{sec:parallel}

The Newton-Kaczmarz method is iterative and, within one iteration, it operates with each individual equation (or record in the case of the application to the Kolmogorov-Arnold model). This means that the standard implementation will not use the capabilities provided by multi-threading that is intrinsic to contemporary computers. The algorithm, however, can be modified such that it becomes suitable for parallel implementation.

In the case of the Kolmogorov-Arnold model with somewhat more complex basis functions, for example, cubic-spline basis functions, most computational resources are spent on the basis function evaluation. Therefore, their parallel evaluation through the entire dataset should be the goal. 

Inner basis functions $\phi_p$ can already be calculated in a parallel way as their arguments are just inputs $x_{i,j}$ that are known \emph{a priori}. This can be done once for the entire dataset and their values can be kept in memory for any consecutive pass of the training algorithm through the dataset.

The problem lies in evaluation of outer basis functions $\psi_l$ as their argument contains parameters $H_{kjp}^q$ that change each iteration. The main idea is to assume that such change is small and to linearise functions $\psi_l$ with respect to parameters $H_{kjp}$ at some value of these parameters $H_{kjp}^\omega$. This decouples the evaluation of $\psi_l$ and the update of $H_{kjp}$, allowing the parallel calculation of the former.

To write the modified algorithm formally, first, equations \eqref{eq:iter0} and \eqref{eq:iter} are rewritten, showing explicitly the parameters and the inputs used for the model evaluation:
\begin{align*}
  &H_{kjp}^{q+1} = H_{kjp}^q - \mu \left(\tilde{y}\left(\vn{H}^q,\vn{G}^q,\vn{x}_i\right) - y_i\right)\zeta\left(\vn{H}^q,\vn{G}^q,\vn{x}_i\right)^{-1} \frac{\partial \tilde{y}\left(\vn{H}^q,\vn{G}^q,\vn{x}_i\right)}{\partial H_{kjp}} , \\
  &G_{kl}^{q+1} = G_{kl}^q - \mu \left(\tilde{y}\left(\vn{H}^q,\vn{G}^q,\vn{x}_i\right) - y_i\right)\zeta\left(\vn{H}^q,\vn{G}^q,\vn{x}_i\right)^{-1} \frac{\partial \tilde{y}\left(\vn{H}^q,\vn{G}^q,\vn{x}_i\right)}{\partial G_{kl}} , 
\end{align*}
where $\vn{H}^q$ and $\vn{G}^q$ stand for the entire sets of parameters $H_{kjp}^q$ and $G_{kl}^q$, respectively. The model output (and its derivatives) is calculated using the parameters at iteration $q$ and the input vector of record $i$, which is denoted as $\tilde{y}\left(\vn{H}^q,\vn{G}^q,\vn{x}_i\right)$. The linearisation of $\tilde{y}$ with respect to the parameters results in
\begin{equation}
\begin{split}
  &\tilde{y}\left(\vn{H}^q,\vn{G}^q,\vn{x}_i\right) = \tilde{y}\left(\vn{H}^\omega,\vn{G}^\omega,\vn{x}_i\right) + \sum_{k,j,p} \frac{\partial \tilde{y}\left(\vn{H}^\omega,\vn{G}^\omega,\vn{x}_i\right)}{\partial H_{kjp}} \left(H_{kjp}^q - H_{kjp}^\omega\right) + \\ 
  &\quad\quad+\sum_{k,l} \frac{\partial \tilde{y}\left(\vn{H}^\omega,\vn{G}^\omega,\vn{x}_i\right)}{\partial G_{kl}} \left(G_{kl}^q - G_{kl}^\omega\vphantom{H_{kjp}^q}\right) .
\end{split}
\end{equation}
Substituting the obtained linearisation into the iterative formulas results in
\begin{align}
  &H_{kjp}^{q+1} = H_{kjp}^q - \mu L_i \zeta\left(\vn{H}^\omega,\vn{G}^\omega,\vn{x}_i\right)^{-1} \frac{\partial \tilde{y}\left(\vn{H}^\omega,\vn{G}^\omega,\vn{x}_i\right)}{\partial H_{kjp}} , \label{eq:parb} \\
  &G_{kl}^{q+1} = G_{kl}^q - \mu L_i \zeta\left(\vn{H}^\omega,\vn{G}^\omega,\vn{x}_i\right)^{-1} \frac{\partial \tilde{y}\left(\vn{H}^\omega,\vn{G}^\omega,\vn{x}_i\right)}{\partial G_{kl}} , \\
  &L_i = \tilde{y}\left(\vn{H}^\omega,\vn{G}^\omega,\vn{x}_i\right) - y_i + \sum_{k,j,p} \frac{\partial \tilde{y}\left(\vn{H}^\omega,\vn{G}^\omega,\vn{x}_i\right)}{\partial H_{kjp}} \left(H_{kjp}^q - H_{kjp}^\omega\right) + \nonumber\\ 
  &\quad\quad+\sum_{k,l} \frac{\partial \tilde{y}\left(\vn{H}^\omega,\vn{G}^\omega,\vn{x}_i\right)}{\partial G_{kl}} \left(G_{kl}^q - G_{kl}^\omega\vphantom{H_{kjp}^q}\right) . \label{eq:pare}
\end{align}
It can be seen that in equations \eqref{eq:parb}-\eqref{eq:pare}, all function evaluations are done using parameters $H_{kjp}^\omega$ and $G_{kl}^\omega$. These parameters can be `frozen' for one pass through the dataset; thus, allowing all function evaluations to be performed in a parallel way. Formally, this means that index $\omega$ changes only once every pass through the dataset: $\omega \in \left\lbrace 0, N, 2N, 3N, \ldots \right\rbrace$.

Iterative formulas \eqref{eq:parb}-\eqref{eq:pare} have been implemented by the authors in the MATLAB version of the code using vectorisation. Thus, computationally-heavy basis function evaluations are out of the main parameter-update loop. Preliminary tests show exactly the same performance in terms of accuracy with the expected significant decrease in computational time.

\subsubsection{Further acceleration of the algorithm}
\label{sec:accel}

The combination of the Kolmogorov-Arnold model with arbitrary basis functions and the Newton-Kaczmarz method for parameter estimation is very flexible, creating an opportunity for accelerating the training method by introducing specific algorithmic enhancements. In particular, this subsection presents three significant enhancements, which were implemented by the authors in their C++ code (but not in the MATLAB code).

\textbf{Gradual increase of the resolution.} The algorithm updates the parameters record-by-record --- if applied straightforwardly, the model's architecture stays constant throughout the entire training. However, it is not necessary to keep it constant. It is possible to start from a coarse model, consisting of the minimal number of the basis functions, and increase the number of the basis functions each time a certain number of iterations is completed. When the number of the basis functions is changed, the parameters must be resampled, but this is a straightforward step that is done for each inner function ($f_{kj}$) and each outer function ($\varPhi_k$) of the model separately. This continues until the predefined number of the basis functions is reached.

Such strategy brings two benefits. It is particularly advantages in terms of resulting accuracy when the inputs of the data records of the dataset do not cover the entire input space, while the basis functions have compact support. In this case, the parameters corresponding to the basis functions, which cover the data gaps, will be obtained via interpolation of coarser models, rather than from random assignment if such strategy is not employed. The second advantage is that the training process is partially done for smaller models, resulting in faster training.

\textbf{Replacement of the basis functions during training.} The models' training time significantly depends on the choice of the basis functions. The simplest functions (i.e. piecewise-linear) require the least amount of calculations, resulting in the fastest training time. As above, it is not necessary to keep the model's architecture constant throughout the entire training. One can perform the majority of the passes through the dataset using the simplest functions (i.e. piecewise-linear) and then switch to more elaborate functions with higher descriptive capabilities (i.e. splines) by reassigning the parameter values directly from one model to another, relying on similarity of the functions. This allows obtaining an accurate model (e.g. with spline basis functions), without using computationally-expensive basis functions throughout the entire training.

\textbf{Skipping of the records.} Since the algorithm iterates record-by-record, it is not necessary to use all records within one pass over the dataset. The parameters' updates are proportional to residuals $\left(\tilde{y}_i - y_i\right)$. For records with small errors $\left|\tilde{y}_i - y_i\right|$, the update steps will not change the parameter values significantly, compared to records with large errors. Therefore, within one pass, the errors can be saved in memory, and within the subsequent passes, only records with errors above the predefined threshold can be updated. This can then be followed one pass where all records are used, recalculating the saved errors and starting another set of passes with records' skipping. Depending on the choice of the threshold, such strategy can lead to significant decrease of the training time.

\section{Data-driven solution of non-linear PDEs}
\label{sec:pdes}

Most multi-physics problems are described by partial differential equations (PDEs). Solution of PDEs provides physical quantities within predefined spatial and temporal domains, but it requires exact boundary data, which is not often readily available due to experimental limitations. Noisy and partially-available data can be utilised in combination with traditional solution methods after post-processing, but such techniques typically have an \emph{ad hoc} nature.

An alternative approach to predicting evolution of physical quantities under scarcity and/or uncertainty of boundary data is physics-informed machine learning \cite{Karniadakis2021}. It is an emerging field focusing on embedding physical constraints (e.g. in form of PDEs) into machine-learning models by training the latter not only on the experimentally-obtained data, but also on the additional data resulting from the enforcement of the constraints. Neural networks are predominantly utilised as the underlying model, and physics-informed neural networks (PINNs) have been successfully used to solve, for example, non-linear PDEs with shock waves \cite{Raissi2019} and the Navier-Stokes equations \cite{Raissi2020}.

Conceptually, the approach consists in assuming that the unknown function of a PDE is the output of a model, the inputs of which are the spatial and the temporal coordinates. The model parameters are then estimated using points internal to the computational domain, at which the residuals are calculated using the PDE, and points at the boundary of the computational domain, at which the residuals are obtained from the data.

\subsection{Kolmogorov-Arnold model as a PDE solver}
\label{sec:KAPDE}

The Kolmogorov-Arnold model in combination with the Newton-Kaczmarz method for parameter estimation can be used not only for data modelling, as described above, but also for data-driven solution of non-linear PDEs. Since the present paper focuses on a model with a single output, the present section is limited to PDEs with scalar unknown functions. The generalisation to PDEs with vector unknown functions is straightforward --- a separate model can be associated with each component of the vector function, which is demonstrated in the next section using a numerical example.

Spatio-temporal domain $\varOmega \in \mathbb{R}^m$ is considered, and $u: \varOmega \to \mathbb{R}$ is the unknown function of partial differential equation
\begin{equation}
  L_\mathrm{I}\left(\vn{x},u,\nabla u,\nabla\nabla u,\ldots\right) = 0 , \quad\quad u=u\left(\vn{x}\right) , \quad\quad \vn{x}\in\varOmega ,
  \label{eq:PDEi}
\end{equation}
where $L_\mathrm{I}$ is some function, $\nabla u$ denotes the vector of partial derivatives $\partial u / \partial x_j$, each component $x_j$ of vector $\vn{x}$ may stand for spatial or temporal coordinate, $\nabla\nabla u$ denotes the matrix of the second partial derivatives, etc. The boundary conditions are written in a similar form:
\begin{equation}
  L_\mathrm{B}\left(\vn{x},u,\nabla u,\nabla\nabla u,\ldots\right) = 0 , \quad\quad \vn{x}\in\varGamma , \quad\quad \varGamma\in\partial\varOmega .
  \label{eq:PDEb}
\end{equation}
It is assumed that functions $L_\mathrm{I}$ and $L_\mathrm{B}$ are such that the problem is well-posed.

A set of randomly-generated domain-internal points $\vn{x}_i\in\varOmega/\varGamma$ for $i\in\left\lbrace 1,\ldots,N_\mathrm{I}\right\rbrace$ is created. In the case of experimentally-measured boundary data, the values of $L_\mathrm{B}$ are available only at a finite number of points, which can be written as $\vn{x}_i\in\varGamma$ for $i\in\left\lbrace N_\mathrm{I}+1,\ldots,N_\mathrm{I}+N_\mathrm{B}\right\rbrace$. The key idea of the approach is to assume that $u$ can be adequately approximated by the Kolmogorov-Arnold model:
\begin{equation}
  u\left(\vn{x}_i\right) = \tilde{y}_i ,
\end{equation}
where $\tilde{y}_i$ is given by equations \eqref{eq:KAU}, \eqref{eq:KAfunc1}, \eqref{eq:KAfunc2}. 

The problem can then be formulated as finding parameters $H_{kjp}$ and $G_{kl}$, such that the absolute values of functions $L_\mathrm{I}$ and $L_\mathrm{B}$ evaluated at points $\vn{x}_i$ are minimised. As above, this problem can be written as minimisation problem
\begin{equation}
  \left\lbrace H_{kjp}^* , G_{kl}^* \right\rbrace = \underset{\vn{H},\vn{G}}{\operatorname{argmin}} \left(\frac{\nu}{N_\mathrm{I}}\sum_{i=1}^{N_\mathrm{I}} L_\mathrm{I}^2\left(\vn{x}_i\right) + \frac{1-\nu}{N_\mathrm{B}}\sum_{i=N_\mathrm{I}+1}^{N_\mathrm{I}+N_\mathrm{B}} L_\mathrm{B}^2\left(\vn{x}_i\right)\right) ,
  \label{eq:KAPDEmin}
\end{equation}
where $H_{kjp}^*$ and $G_{kl}^*$ are the sought parameters, $\nu$ is a weight factor associated with the internal data.

An approximate solution of problem \eqref{eq:KAPDEmin} can be found using the Newton-Kaczmarz method. Equations \eqref{eq:PDEi} and \eqref{eq:PDEb} evaluated at points $\vn{x}_i$ are rewritten with respect to the model parameters:
\begin{align}
  &L_\mathrm{I}\left(\vn{x}_i\right) = L_\mathrm{I}\left(\vn{H}, \vn{G}, \vn{x}_i\right) = 0 , 
    \quad\quad i\in\left\lbrace 1,\ldots,N_\mathrm{I}\right\rbrace , \label{eq:PDEip} \\
  &L_\mathrm{B}\left(\vn{x}_i\right) = L_\mathrm{B}\left(\vn{H}, \vn{G}, \vn{x}_i\right) = 0 ,
    \quad\quad i\in\left\lbrace N_\mathrm{I}+1,\ldots,N_\mathrm{I}+N_\mathrm{B}\right\rbrace .
\label{eq:PDEbp}
\end{align}
where $\vn{H}$ and $\vn{G}$ stand for the entire sets of parameters $H_{kjp}$ and $G_{kl}$, respectively. This is a system of non-linear equations. As above, the solution process starts with random initial values $H_{kjp}^0$ and $G_{kl}^0$ and updates the parameters iteratively point-by-point in the opposite direction to the gradient of $L$ with respect to the parameters. The subscripts of $L$ are omitted for brevity --- it is implied that $L = L_\mathrm{I}\left(\vn{x}_i\right)$ when $\vn{x}_i\in\varOmega/\varGamma$ and $L = L_\mathrm{B}\left(\vn{x}_i\right)$ when $\vn{x}_i\in\varGamma$. This gives the following iterative procedure:
\begin{equation}
  H_{kjp}^{q+1} = H_{kjp}^q - \mu L \zeta^{-1} \frac{\partial L}{\partial H_{kjp}} , \quad\quad\quad
  G_{kl}^{q+1} = G_{kl}^q - \mu L \zeta^{-1} \frac{\partial L}{\partial G_{kl}} , \label{eq:iterPDE}
\end{equation}
where $\mu\in\left(0,1\right]$ is the numerical damping parameter. Superscript $q$ denotes the approximation of the parameters at iteration $q$. The values of $L$ and its derivatives are calculated using the parameters at iteration $q$. Variable $\zeta$ is the squared norm of the gradient of $L$ with respect to the parameters:  
\begin{equation}
  \zeta = \sum_{k=1}^{d} \sum_{j=1}^{m} \sum_{p=1}^{n} \left(\frac{\partial L}{\partial H_{kjp}}\right)^2 + \sum_{k=1}^{d} \sum_{l=1}^{s} \left(\frac{\partial L}{\partial G_{kl}}\right)^2 .
\end{equation}
As previously, index $i$ changes each iteration. 

In the case of data-driven solution of PDEs, the set of equations for parameter estimation consists of two distinct groups, corresponding to the internal and to the boundary points; therefore, the choice of how index $i$ changes will influence the result, i.e. the user might decide to pick the internal or the boundary points more frequently. Weight factor $\nu$ from problem \eqref{eq:KAPDEmin} corresponds to the frequency of selection of a record with domain-internal input $\vn{x}_i$. To perform the model training systematically, consecutive iterations of the model update can be combined into batches. One such batch update consists of using $N_\mathrm{b}$ boundary points and $N_\mathrm{i}$ internal points. In this case, $\nu = N_\mathrm{i}/\left(N_\mathrm{i}+N_\mathrm{b}\right)$.

\section{Numerical examples}
\label{sec:simulation}

In this section, several numerical examples are presented to show the efficiency of the proposed approach. After having represented unknown functions $f_{kj}$ and $\varPhi_k$ via basis functions $\phi_p$ and $\psi_l$ and parameters $H_{kjp}$ and $G_{kl}$, various methods can be used to determine the parameters. The straightforward choice is the Gauss-Newton (GN) method. Therefore, it is necessary to investigate how it performs compared to the Newton-Kaczmarz (NK) method for the considered problem, which is done in the first example. Parameter estimation of the Kolmogorov-Arnold model using the GN method is summarised in appendix \ref{sec:app}. Furthermore, as the NK method is related to the stochastic gradient descent (SGD) method (as summarised in appendix \ref{sec:SGD}), it is useful to compare the NK method to other variants of the SGD method. For such comparison, the relatively-popular `Adam' SGD method is chosen \cite{Kingma2015}. For completeness of the paper, it is summarised in appendix \ref{sec:adam}. In the examples below, the non-parallel version of the NK method is used.

\subsection{Example 1: Parameter estimation of a ridge function}
\label{sec:examp1}

To compare the parameter estimation methods, the ridge function model, equation \eqref{eq:ridge}, is taken. It can be viewed as a particular case of the Kolmogorov-Arnold model with one addend, $d=1$, and with the inner functions being linear, $n=1$ and $\phi_p\left(x\right) = x$. Such model does not have redundancy, and the parameter estimation problem has one well-defined global minimum. It can also have local minima that can `trap' the parameter estimation method. The idea behind this example is to set up the most challenging case, such that the outer function is highly non-liner and convergence to local minima results in insufficient accuracy. It should be noted of course that in complex real-world scenarios, getting to some local minima may perfectly satisfy the user, but it is useful to know that some training methods are better than others in terms of getting to the global minimum. In particular, since Newton-type methods require an initial guess, it is important that the model training technique is robust with respect to the choice of the guess and does not get stuck at local minima. Therefore, this example aims to compare the dependence of the methods' convergence on the initial guess.

The example consists in assuming a certain model, generating the data (i.e. generating the inputs and calculating the corresponding outputs using the assumed model), assuming that the model parameters are now unknown, and estimating them using the data. The advantage of such approach is the availability of the exact solution of the parameter estimation problem. Of course, the structure of the new model (e.g. the number of the basis functions) is chosen to be identical to the original model that generates the data because the objective is to recover the original model. The code of this example is available via GitHub\footnote{https://github.com/andrewpolar/RidgeIdentM}.

As before, outer function $\varPhi$ is represented via the basis functions and the parameters --- equation \eqref{eq:KAfunc2} with omitted index $k$. The Gaussian basis functions, equation \eqref{eq:gaussbf}, with $\gamma=2$, $t_\mathrm{min}=0.5$, and $t_\mathrm{max} = 2.5$, are taken. To generate the data, a model with $m=5$ inputs and $s=3$ basis functions is taken. The following parameters are assumed:
\begin{equation}
  c = \begin{bmatrix} -0.7 & 2.5 & -1.2 & 0.8 & 1.6 \end{bmatrix}, \quad\quad 
  G = \begin{bmatrix} 2.1 & -0.9 & 0.7 \end{bmatrix} .
  \label{eq:model}
\end{equation}
To create the datasets of $N=400$ records, uniformly-distributed random inputs from $x_\mathrm{min} = 0$ to $x_\mathrm{max} = 1$ are generated, $x_{i,j} \sim \operatorname{unif}\left(0,1\right)$, and outputs $y_i$ are calculated. 

Having the datasets, the parameter estimation is done using the NK, the GN, and the `Adam' SGD methods. In the NK method, numerical damping parameter $\mu$ is varied (the values are given in the corresponding tables) and $10^4$ iterations are performed. In the GN method, tolerance $\delta = 10^{-12}$ is set and the maximum number of iterations is taken to be $100$; furthermore, the initial value of $\mu$ is set to $0.1$, as soon as $\| \vn{Z}^{q+1} - \vn{Z}^q \| < 0.1$ is fulfilled, the value of $\mu$ is changed to $1$; moreover, to obtain the most accurate Hessian matrix and to maximise the performance of the method, the second derivatives of $\tilde{y}$ are not neglected; finally, in this example, the regularisation is not used, $\lambda = 0$. In the `Adam' SGD method, $10^4$ iterations are performed (the same as in the NK method to establish the proper comparison), parameter $\varepsilon = 10^{-8}$ is set (the default value proposed in \cite{Kingma2015}), and parameters $\beta_1,\beta_2,\mu$ are varied (the values are given in the corresponding tables).
 
The initial guess for the parameters is taken as the exact solution plus a random perturbation, i.e.
\begin{equation}
  c_j^0 = c_j + \alpha a_{j} , \quad\quad 
  G_l^0 = G_l + \alpha b_{l} , 
\end{equation}
where $a_{j},b_{l} \sim \operatorname{unif}\left(-0.5,0.5\right)$ are uniformly-distributed random numbers. For small $\alpha$, the initial guess is close to the true solution and the methods are expected to converge, while for large $\alpha$, the methods might not converge.

The normalised root-mean-square error (RMSE),
\begin{equation}
  E_\mathrm{RMSE} = \frac{1}{y_\mathrm{max}-y_\mathrm{min}}\sqrt{ \frac{1}{N}\sum_{i=1}^N \left(y_i - \tilde{y}_i\right)^2 } ,
  \label{RMSE}
\end{equation}
is used as an accuracy metric. After the methods produce the model parameters for given $\alpha$, error $E_\mathrm{RMSE}$ is calculated. Since the generation of the datasets and the initial guess for the parameters involve randomisation, error $E_\mathrm{RMSE}$ can significantly differ from one run\footnote{One execution of the code is implied by the `run'.} to another; therefore, its average behaviour should be considered. In particular, an illustrative characteristic of the performance of the methods is the average proportion of runs, in which $E_\mathrm{RMSE}$ is below some threshold. Such characteristic is more useful than a direct comparison of average $E_\mathrm{RMSE}$ produced by the GN and the NK methods because it is not sensitive to the outliers, where term `outlier' is used loosely to refer to the results way above the threshold.

The average percentages of runs, for which the RMSE is below $5\%$, $10\%$, and $20\%$ for different amplitude of the perturbations of the initial guess, parameter $\alpha$, are shown in table \ref{tab:GNres} for all three methods. As expected, as $\alpha$ increases, the percentage of the accurate runs drops. More importantly, it can be seen that the GN method consistently results in fewer runs with the error below a given threshold than the NK method. Furthermore, the accuracy of the GN method drops more rapidly than of the NK method. For example, for $\alpha = 1.2$, more than $90\%$ of the runs of the NK method produce errors less than $20\%$, while only about $50\%$ of the runs of the GN method have such accuracy; for $\alpha = 2.8$, which is close to a completely random initial guess, since such value of $\alpha$ is greater than any parameter value in the assumed model, about $20\%$ of the runs of the NK method produce errors less than $10\%$, while only about $1\%$ of the runs of the GN method have such accuracy.

\begin{table}
  \begin{center}
    \begin{tabular}{|l|lllllll|}
    \hline
    $\alpha$ & $0.4$ & $0.8$ & $1.2$ & $1.6$ & $2.0$ & $2.4$ & $2.8$ \\
    \hline
    \multicolumn{8}{|c|}{Gauss-Newton (GN) method} \\
    \hline
    err. $<5\%$ & $93.6 \pm 1.3$ & $59.0 \pm 5.2$ & $31.4 \pm 1.7$ & $12.4 \pm 2.1$ & $5.0 \pm 2.1$ & $2.0 \pm 1.6$ & $0.4 \pm 0.5$ \\
    err. $<10\%$ & $93.8 \pm 1.1$ & $59.0 \pm 5.2$ & $32.6 \pm 2.3$ & $13.6 \pm 2.3$ & $5.6 \pm 1.9$ & $2.2 \pm 1.9$ & $1.2 \pm 0.8$ \\
    err. $<20\%$ & $99.2 \pm 0.4$ & $80.0 \pm 4.8$ & $55.0 \pm 6.4$ & $33.8 \pm 2.9$ & $21.8 \pm 4.8$ & $11.0 \pm 3.1$ & $7.4 \pm 2.9$ \\
    \hline
    \% of conv. & $99.4 \pm 0.5$ & $83.2 \pm 5.5$ & $62.4 \pm 7.2$ & $40.8 \pm 3.6$ & $29.8 \pm 7.2$ & $22.2 \pm 2.9$ & $19.6 \pm 5.6$ \\
    RMSE in \% & $0.8 \pm 0.2$ & $4.3 \pm 0.1$ & $8.3 \pm 1.1$ & $11.9 \pm 1.0$ & $15.1 \pm 2.8$ & $22.1 \pm 2.5$ & $24.7 \pm 3.1$ \\
    \hline
    \multicolumn{8}{|c|}{Newton-Kaczmarz (NK) method, $\mu = 0.1$} \\
    \hline
    err. $<5\%$ & $98.4 \pm 0.9$ & $78.2 \pm 3.0$ & $49.4 \pm 4.7$ & $35.8 \pm 4.1$ & $24.8 \pm 2.6$ & $16.0 \pm 5.1$ & $13.0 \pm 4.3$ \\
    err. $<10\%$ & $99.8 \pm 0.4$ & $95.0 \pm 2.0$ & $78.0 \pm 7.2$ & $56.2 \pm 3.6$ & $42.2 \pm 4.3$ & $30.2 \pm 5.8$ & $20.4 \pm 4.5$ \\
    err. $<20\%$ & $100 \pm 0.0$ & $98.8 \pm 1.1$ & $90.8 \pm 3.8$ & $75.4 \pm 1.7$ & $64.4 \pm 3.2$ & $48.0 \pm 5.3$ & $34.8 \pm 7.2$ \\
    \hline
    \multicolumn{8}{|c|}{Newton-Kaczmarz (NK) method, $\mu = 1$} \\
    \hline
    err. $<5\%$ & $96.6 \pm 2.3$ & $82.0 \pm 6.0$ & $67.6 \pm 2.7$ & $50.6 \pm 8.4$ & $37.0 \pm 3.9$ & $23.8 \pm 4.9$ & $17.6 \pm 3.0$ \\
    \hline
    \multicolumn{8}{|c|}{`Adam' SGD method, $\mu = 10^{-3}$, $\beta_1 = 0.9$, $\beta_2 = 0.999$ (default parameters)} \\
    \hline
    err. $<5\%$ & $98.6 \pm 1.3$ & $79.6 \pm 5.3$ & $51.6 \pm 8.1$ & $28.0 \pm 3.1$ & $15.6 \pm 3.2$ & $6.2 \pm 0.8$ & $1.8 \pm 1.3$ \\
    \hline
    \multicolumn{8}{|c|}{`Adam' SGD method, $\mu = 10^{-2.5}$, $\beta_1 = 0.9$, $\beta_2 = 0.99$ (fine-tuned parameters)} \\
    \hline
    err. $<5\%$ & $96.2 \pm 1.3$ & $85.6 \pm 4.6$ & $72.8 \pm 5.8$ & $49.4 \pm 3.5$ & $37.8 \pm 2.9$ & $26.6 \pm 3.5$ & $14.0 \pm 1.9$ \\
    \hline
    \end{tabular}
  \end{center}
  \caption{The average number of runs (executions of the code) out of $100$, where the normalised root-mean-square error is below $5\%$, $10\%$, $20\%$ for different amplitude of the perturbations of the initial guess (parameter $\alpha$), for the GN method (lines 3-5), for the NK method with $\mu = 0.1$ (lines 9-11), for the NK method with $\mu = 1$ (line 13), for the `Adam' SGD method with default parameters (line 15) and for the `Adam' SGD method with fine-tuned parameters (line 17). The average number of runs out of $100$, where the convergence criteria of the GN method is fulfilled for different $\alpha$ (line 6). The average normalised root-mean-square error in $\%$, calculated only for the converged runs of the GN method, for different $\alpha$ (line 7). In each cell, the mean and the standard deviation across $5$ experiments, each consisting of $100$ runs, are given.}
  \label{tab:GNres}
\end{table}

The percentages of the accurate runs of the methods might seem to be relatively small, but this example is specifically chosen to be challenging for the methods. In particular, the exact values of $G_l$ given in equation \eqref{eq:model} are such that outer function $\varPhi$ is highly non-linear and non-monotonous. Furthermore, since the model is the ridge function, there is no redundancy --- it consists of only one addend, parameters of which must be obtained precisely, otherwise the error is large. This means that the convergence of the methods to local minima (with the obtained parameters that are different from the exact model parameters) is insufficient in terms of accuracy.  

In the comparison above, all runs of the GN method are used, even if the convergence criteria is not fulfilled and the iterations are stopped after the maximum number of iterations has been reached. An alternative approach is disregarding the runs where the convergence criteria is not fulfilled. However, this implies that if the convergence criteria is fulfilled, then the result should be considered to be accurate. In this case, first, the percentage of runs with the `converged' results drops with $\alpha$, as shown in line 6 of table \ref{tab:GNres}. Second, the average RMSE, calculated only for the `converged' runs, still increases with $\alpha$ and reaches almost $25\%$ for $\alpha = 2.8$, as shown in line 7 of table \ref{tab:GNres}. This means that the runs converge to local minima, providing insufficiently accurate result.

Table \ref{tab:GNres} also provides a comparison between the NK and the `Adam' SGD methods. The latter has three fine-tuning parameters $\beta_1,\beta_2,\mu$, and it is possible to adjust those such that the method performs close to the NK method in this particular problem. The results with the default parameters (as proposed in \cite{Kingma2015}) and the results with the fine-tuned parameters are provided in table \ref{tab:GNres}. The fine-tuning has been done by the authors `knowing' the exact solution\footnote{The numerical parameters of the method were varied until the best solution was achieved. Because the data has been generated by the originally-known model, when the `Adam' SGD method is used to train another model, it is known exactly what result should be obtained.}; thus, it is only provided to show that the optimal values of the parameters exist. The performance of the `Adam' SGD method with the default parameters is close to the performance of the NK method with $\mu=0.1$ (although the NK method is noticeably better when $\alpha$ is large). The `Adam' SGD method with the fine-tuned parameters is close to the NK method with $\mu=1$. The small difference is not surprising given the relationship between the NK and the SGD methods.

There is significant sensitivity of the performance of the `Adam' SGD method to the value of the learning rate. In particular, in table \ref{tab:Adam}, the same accuracy metric (the average percentages of runs where the RMSE is below $5\%$) is compared across the methods with $\mu$ varied. For the considered problem, the `Adam' SGD method provides good results only when the learning rate is between $10^{-2}$ and $10^{-3}$. For $\mu = 10^{-3.5}$, there is a significant drop in the number of accurate runs when $\alpha = 1.6$, i.e. when the initial guess for the solution is relatively far from the exact solution. In this case, the performance of the method is even worse than that of the GN method (when compared to line 3 of table \ref{tab:GNres}). For $\mu = 10^{-1.5}$, there is a very low percentage of accurate runs when $\alpha = 0.4$, i.e. when the initial guess for the solution is very close to the exact solution. For `good' initial guesses, the Newton-type methods almost always converge well, while the `Adam' SGD method with $\mu = 10^{-1.5}$ probably overshoots the optimal solution update due to the large step size. 

\begin{table}
  \begin{center}
    \begin{tabular}{|l|lllll|}
    \hline
    \multicolumn{6}{|c|}{Newton-Kaczmarz (NK) method} \\
    \hline
    $\mu$ & $10^0$ & $10^{-0.5}$ & $10^{-1}$ & $10^{-1.5}$ & $10^{-2}$ \\
    \hline
    err. $<5\%$, $\alpha = 0.4$ & $96.6 \pm 2.3$ & $98.4 \pm 0.9$ & $98.4 \pm 0.9$ & $97.0 \pm 2.2$ & $92.6 \pm 1.3$ \\
    err. $<5\%$, $\alpha = 1.6$ & $50.6 \pm 8.4$ & $46.4 \pm 6.5$ & $35.8 \pm 4.1$ & $24.2 \pm 1.3$ & $9.4 \pm 1.8$ \\
    \hline
    \multicolumn{6}{|c|}{`Adam' SGD method, $\beta_1 = 0.9$, $\beta_2 = 0.999$} \\
    \hline
    $\mu$ & $10^{-1.5}$ & $10^{-2}$ & $10^{-2.5}$ & $10^{-3}$ & $10^{-3.5}$ \\
    \hline
    err. $<5\%$, $\alpha = 0.4$ & $31.2 \pm 3.6$ & $75.6 \pm 3.0$ & $97.0 \pm 1.2$ & $98.6 \pm 1.3$ & $95.6 \pm 3.0$ \\
    err. $<5\%$, $\alpha = 1.6$ & $21.0 \pm 3.7$ & $45.6 \pm 3.9$ & $44.8 \pm 7.0$ & $28.0 \pm 3.1$ & $7.6 \pm 2.7$ \\
    \hline
    \end{tabular}
  \end{center}
  \caption{The average number of runs (executions of the code) out of $100$, where the normalised root-mean-square error is below $5\%$ for different amplitude of the perturbations of the initial guess, parameter $\alpha$, and for different values of parameter $\mu$ for the NK method (lines 3-4) and for the `Adam' SGD method (lines 7-8). In each cell, the mean and the standard deviation across $5$ experiments, each consisting of $100$ runs, are given.}
  \label{tab:Adam}
\end{table}

The NK method has only one fine-tuning parameter --- $\mu$. In the `Adam' SGD method it is referred to as the learning rate, while in the NK method it is more logical to call it the numerical damping parameter (as done in the preceding sections). There is an intuitive understanding that for the exact data (i.e. the input-output data corresponding to the consistent system of non-linear equations that has a unique solution), the optimal value of $\mu$ in the NK method is $1$. This follows from the error analysis (appendix \ref{sec:errAn}), where it is shown that such value of $\mu$ corresponds to the fastest error decay. As the example of this section deals with the exact data, $\mu = 1$ straightforwardly gives the best performance of the NK method, as seen in tables \ref{tab:GNres} and \ref{tab:Adam}. Using $\mu < 1$ makes the method robust for the case of noisy data (some insights how this works are provided in \cite{Poluektov2020}, where such numerical damping has been studied within the context of the Urysohn model with the Kaczmarz-based model training algorithm).

The percentages in tables \ref{tab:GNres} and \ref{tab:Adam} are calculated based on $5$ separate experiments, each consisting of $100$ runs. Within each experiment, the number of runs with the error below a given threshold is calculated; the mean and the standard deviation across the experiments are given in the tables. The same is done for the percentage of the converged runs of the GN method and the corresponding RMSEs.

\subsection{Example 2: Approximation of a non-linear function}

The purpose of this example is to demonstrate that the NK method can indeed solve the parameter estimation problem for the Kolmogorov-Arnold model approximating a non-linear multivariate function. The code of this example is available via GitHub\footnote{https://github.com/andrewpolar/UStandardM}. The `full' model is taken with the piecewise-linear basis functions.

The following non-linear function of $m=5$ inputs is chosen:
\begin{equation}
\begin{split}
  &y_i = \frac{2 + 2 x_{i,3}}{3\pi} \left( \operatorname{arctan}\left( 20\left( x_{i,1} - \frac{1}{2} + \frac{x_{i,2}}{6}\right)\exp\left(x_{i,5}\right) \right) + \frac{\pi}{2} \right) + \\ 
  &\vphantom{y} + \frac{2 + 2 x_{i,4}}{3\pi} \left( \operatorname{arctan}\left( 20\left( x_{i,1} - \frac{1}{2} - \frac{x_{i,2}}{6}\right)\exp\left(x_{i,5}\right) \right) + \frac{\pi}{2} \right) .
\end{split}
\label{eq:formula2}
\end{equation}
The datasets consisting of $N = 10^4$ records for the training and of $N_\mathrm{val} = 2000$ records for the validation are generated. This involves taking uniformly-distributed random inputs $x_{i,j} \sim \operatorname{unif}\left(0,1\right)$ and calculating outputs $y_i$ using equation \eqref{eq:formula2}. In this example, the models with the piecewise-linear basis functions are assessed. The numbers of the basis functions were chosen by variation and comparison of the performance on the training dataset, which resulted in $n=5$ inner and $s=7$ outer basis functions for this problem and dataset size.

\begin{figure}
  \begin{center}
    \includegraphics[scale=0.9]{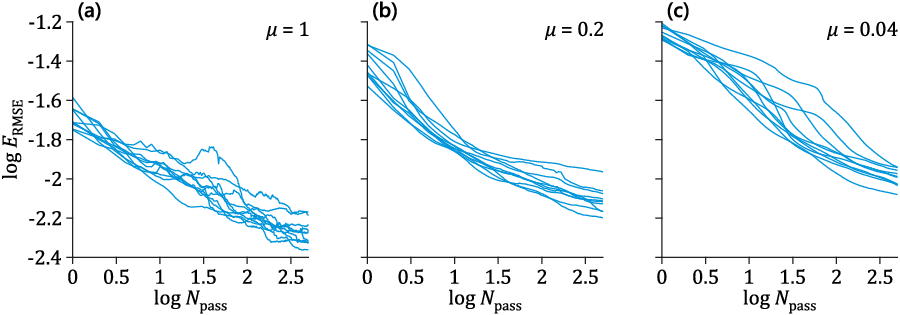}
  \end{center}
  \caption{The dependence of the RMSE on the number of passes through the data for different values of $\mu$.}
  \label{fig:nlFunc}
\end{figure}

Since the method is iterative, it is interesting to study the convergence rate with respect to the number of iterations. The number of iterations should be large, but each iteration is computationally cheap, as each iteration of the NK method operates only with one record. Therefore, the algorithm goes through all the training records consecutively multiple times. After each pass through the training data, the RMSE is calculated using the validation data. Thus, the RMSE can be plotted as a function of the number of passes through the data, which is illustrated in figure \ref{fig:nlFunc}. Three different values of the numerical damping parameter are used; for each case, $10$ runs with different data and with different initial guesses for the parameters are performed.

It can be seen that for any $\mu$, the method converges approximately logarithmically with respect to the number of passes, i.e. there is approximately linear dependence between $\log E_\mathrm{RMSE}$ and $\log N_\mathrm{pass}$. However, smaller $\mu$ leads to higher absolute value of $E_\mathrm{RMSE}$ for the same number of passes. The practical purpose of decreasing $\mu$ is to filter out the noise, which is absent in this example. A more detailed study of the noise filtering abilities of the Kaczmarz method has been performed for the case of the Urysohn model in \cite{Poluektov2020}. Finally, it should be mentioned that the model obtained using the proposed technique is accurate, with the RMSE of around $0.5\%$ for $\mu=1$ after $500$ passes, as seen in figure \ref{fig:nlFunc}.

\subsection{Example 3: Solution of a PDE}

The purpose of this example is to demonstrate that the Kolmogorov-Arnold model with the NK method for parameter estimation can indeed be used for data-driven solution of PDEs. The `full' model is taken with the spline basis functions.

To test the approach properly, it is necessary to consider a PDE with an analytical solution. The following function is taken:
\begin{equation}
  u\left(x_1,x_2\right) = x_1 \exp\left(x_2 - x_2^2 \right) .
\end{equation}
It is easy to construct a second-order PDE, solution of which is function $u$ above:
\begin{equation}
  \frac{\partial^2 u}{\partial x_1^2} + 2 x_1 x_2 \frac{\partial^2 u}{\partial x_1 \partial x_2} + \frac{\partial^2 u}{\partial x_2^2} + 2 x_1 \frac{\partial u}{\partial x_1} - \frac{\partial u}{\partial x_2} = 0 .
  \label{eq:pde2}
\end{equation}
Domain $\left[0,2\right]\times\left[0,2\right]$ is considered, resulting in the boundary conditions for the PDE:
\begin{equation*}
  \left.u\right|_{x_1=0} = 0 , \quad\quad 
  \left.u\right|_{x_1=2} = 2\exp\left(x_2 - x_2^2 \right) , \quad\quad 
  \left.u\right|_{x_2=0} = x_1 , \quad\quad 
  \left.u\right|_{x_2=2} = x_1\exp\left(-2\right) .
\end{equation*}
The structure of the Kolmogorov-Arnold model is chosen to consist of $n=7$ inner and $s=7$ outer spline basis functions. Such structure was obtained by training the model on a dataset of a given size and selecting $n$ and $s$ that resulted in the minimum error compared to the exact solution. Numerical damping parameter $\mu=0.2$ is used. Batch model updates consisted of $N_\mathrm{b} = 4$ boundary points (one point per boundary condition) and $N_\mathrm{i} = 20$ internal points. Thus, weight factor $\nu$ is shifted towards the internal data. 

\begin{figure}
  \begin{center}
    \includegraphics[scale=0.9]{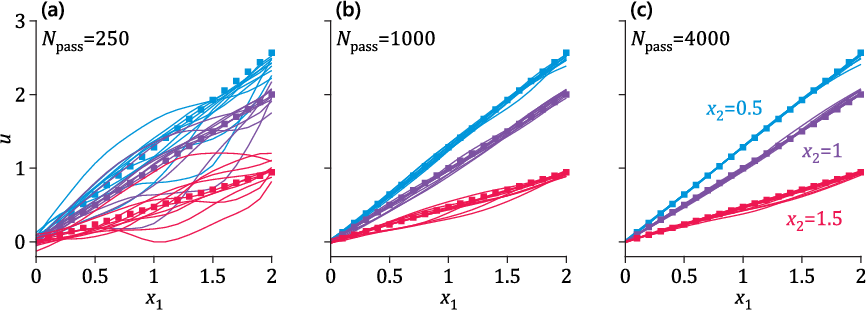}
  \end{center}
  \caption{The sections of function $u\left(x_1,x_2\right)$ at $x_2=0.5$ (blue), $x_2=1$ (brown), and $x_2=1.5$ (red). Solid lines are the numerical solutions of the second-order PDE, which are obtained by training the Kolmogorov-Arnold model using the different number of passes (batch updates). Symbols show the analytical solution.}
  \label{fig:pde}
\end{figure}

The model training process starts with random parameters and converges as the number of batch updates increases. To illustrate this, the number of batch updates is varied and, in each case, $10$ models are trained. Figure \ref{fig:pde} shows the comparison between the analytical solution and the numerical solutions given by the Kolmogorov-Arnold model. It can be seen that the spread of the curves corresponding to the numerical solutions logarithmically decreases with the increase of the number of batch updates. The numerical solutions converge to the analytical solution --- after $10^5$ batch updates, the RMSE of $0.098\% \pm 0.021\%$ is obtained (the numbers are the mean and the standard deviation across $5$ runs). The latter is calculated according to equation \eqref{RMSE} on the structured grid with the spatial step of $0.1$.

\subsection{Example 4: Solid mechanics problem}

The purpose of this example is to expand on the results of the previous example and to show more application-oriented PDE. In particular, the linear momentum balance equation from solid mechanics is considered. This vector equation describes displacement of material points of a deformable solid, which is taken to be isotropic linear elastic:
\begin{equation*}
  \frac{\partial \sigma_{ij}}{\partial x_i} = 0 , \quad\quad
  \sigma_{ij} = \bar{\lambda} \delta_{ij} \varepsilon_{kk} + 2 \bar{\mu} \varepsilon_{ij} , \quad\quad
  \varepsilon_{ij} = \frac{1}{2}\left(\frac{\partial u_i}{\partial x_j} + \frac{\partial u_j}{\partial x_i}\right),
\end{equation*}
where $x_i$ are the spatial coordinates, $u_i$ are the components of the unknown displacement vector, $\bar{\lambda}$ and $\bar{\mu}$ are the Lam\'{e} parameters, $\delta_{ij}$ is the Kronecker delta, and the Einstein summation convention is adopted. In 2D, there are two components of the displacement vector and two equations.

It is easy to verify that solution in the following form satisfies the equation:
\begin{align*}
  &u_1 = \left( \left( C_1 - C_2 \beta \right) e^{\kappa x_1} + C_2 x_1 e^{\kappa x_1} - \left( C_3 + C_4 \beta \right) e^{-\kappa x_1} - C_4 x_1 e^{-\kappa x_1} \right) \sin\left(\kappa x_2\right) , \\
  &u_2 = \left( C_1 e^{\kappa x_1} + C_2 x_1 e^{\kappa x_1} + C_3 e^{-\kappa x_1} + C_4 x_1 e^{-\kappa x_1} \right) \cos\left(\kappa x_2\right) , \quad 
  \beta = \left(\bar{\lambda}+3\bar{\mu}\right)/\left(\kappa\bar{\lambda}+\kappa\bar{\mu}\right) ,
\end{align*}
where $C_1$, $C_2$, $C_3$, $C_4$, $\kappa$ are some constants.

For this example, domain $\left[-1,1\right]\times\left[-1,1\right]$ is considered, $\bar{\lambda}=3$ and $\bar{\mu}=2$ are taken, and to construct the analytical solution, $C_1=-0.2$, $C_2=0.1$, $C_3=-0.07$, $C_4=0$, $\kappa=2.5$ are taken. The displacement boundary conditions are imposed via the expression for the analytical solution.

Two `full' Kolmogorov-Arnold models with the spline basis functions are fit to the data. The first and the second models represent $u_1$ and $u_2$, respectively. To calculate the non-normalised RMSE, a separate structured rectangular grid of points with the step size of $0.01$ is created, and the error is obtained from the difference between the numerical and the analytical solutions.

In numerical analysis, an important property of the method is the convergence rate. The Kolmogorov-Arnold model uses the basis functions with compact support. The latter can also be found in traditional methods, such as the finite-element method (FEM). Hence, one can say that the input coordinates are discretised on structured grids. However, in FEM, a 2D problem is discretised using a 2D mesh. Here, there several 1D meshes, for the input and for the intermediate variables. 

\begin{figure}
  \begin{center}
    \includegraphics[scale=0.9]{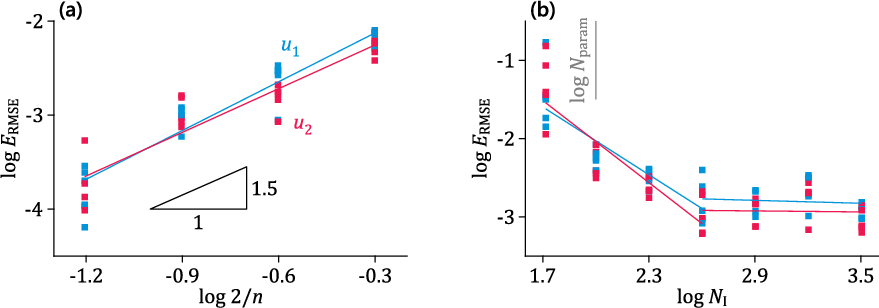}
  \end{center}
  \caption{The dependence of the non-normalised RMSE on the mesh size corresponding to the inputs (a) and on the number of domain-internal points for the training (b). The mesh size is defined as $2/n$, where $2$ is the size of the domain and $n$ is the number of the inner basis functions. The vertical line in (b) indicates $N_\mathrm{I}$ equal to the number of the model's parameters. The red and the blue colours correspond to the solutions for $u_1$ and $u_2$, respectively.}
  \label{fig:mech}
\end{figure}

The convergence can be measured with respect to the number of the basis functions corresponding to the inner and the outer functions. To estimate the former, the models consisting of $n=4,8,16,32$ inner and $s=4$ outer spline basis functions are trained. To train the models, $N_\mathrm{I} = 1400$ domain-internal and $N_\mathrm{B} = 1400$ boundary points are used, and $200$ passes through the data are made. In figure \ref{fig:mech}a, it can be seen that the non-normalised RMSE decays with the rate of $1.5$ with the decrease of the mesh size that corresponds to the input coordinates (i.e. inversely proportional to the number of the inner basis functions). The authors also tried to estimate the convergence rate with respect to the number of the outer basis functions, but the results were inconclusive, requiring more detailed and comprehensive study, which will be the subject of future research.

An interesting question is how many points are required to fit the models. To estimate the latter, the models consisting of $n=8$ inner and $s=4$ outer spline basis functions are trained (i.e. $100$ parameters per model), and the number of domain-internal points is varied, keeping $N_\mathrm{I} = N_\mathrm{B}$. In figure \ref{fig:mech}b, it can be seen that the RMSE first decays with the increase of the number of points, but then stays constant. In this particular example, after the number of domain-internal points reaches four times the number of model's parameters, the accuracy does not change anymore (one domain-internal point corresponds to two equations, but also two models are trained). This ratio may change for other examples, but it can be expected that the trend will be qualitatively the same.

\subsection{Additional examples: Execution times of different codes/methods}

The purpose of this subsection is to compare the execution times of various methods and codes. Of course, the performance of the methods is very much dependent on the specific problem and on the aspects of the written code; therefore, any analysis, including that of the present subsection, will provide only a partial picture, and the codes are not guaranteed to perform the same way in other examples.

\subsubsection{Benchmark 1}

As mentioned in the introduction, since the publication of the first preprint of the present paper, an independent Python implementation of the Kolmogorov-Arnold model/network with the LBFGS (limited-memory Broyden-Fletcher-Goldfarb-Shanno) training method has been posted on GitHub --- the code supplementing reference \cite{Liu2024}. Therefore, the first example of \cite{Liu2024} is taken. It consists in fitting the following function of $m=2$ inputs:
\begin{equation}
  y_i = \exp \left( \sin\left( \pi x_{i,1} \right) + x_{i,2}^2 \right) .
\label{eq:formula3}
\end{equation}
The datasets consisting of $N = 1000$ records for the training and of $N_\mathrm{val} = 1000$ records for the validation are generated. This involves taking uniformly-distributed random inputs $x_{i,j} \sim \operatorname{unif}\left(-1,1\right)$ and calculating outputs $y_i$ using equation \eqref{eq:formula3}. Using the training datasets, the `full' Kolmogorov-Arnold models ($5$ addends) with $n=6$ inner and $s=6$ outer cubic-spline basis functions are constructed using different codes as described below. The summary of the results is provided in table \ref{tab:times1}.

\begin{table}
  \begin{center}
    \begin{tabular}{|l|ll|}
    \hline
    & accuracy ($\%$) & execution time (s) \\
    \hline
    Python code of ref. \cite{Liu2024} & $0.73 \pm 0.31$ & $9.27$ \\
    present paper, MATLAB code, GN method & $0.25 \pm 0.08$ & $0.22$ \\
    present paper, C\# code, NK method & $0.71 \pm 0.13$ & $0.66$ \\
    MATLAB's neural network & $0.66 \pm 0.16$ & $0.47$ \\
    \hline
    \end{tabular}
  \end{center}
  \caption{The average non-normalised RMSEs and the average execution times for different methods and implementations --- benchmark 1. For the accuracy, the mean and the standard deviation across $10$ runs, are given.}
  \label{tab:times1}
\end{table}

The structure of the Kolmogorov-Arnold model (the number of addends and basis functions) follows the setup of the tutorial\footnote{https://kindxiaoming.github.io/pykan/intro.html} that supplements reference \cite{Liu2024}. The only addition is a random seed to see the results changing depending on the initial guess. Training procedure \verb|model.train| is executed specifying $20$ steps of the method.

To facilitate the comparison of the results to other papers, the accuracy given in table \ref{tab:times1} is the non-normalised RMSE or $\left(y_\mathrm{max}-y_\mathrm{min}\right) E_\mathrm{RMSE}$, where $E_\mathrm{RMSE}$ is given by equation \eqref{RMSE}. The values are given in $\%$. All codes were executed on a laptop with Intel Core i5-1145G7 processor.

The first comparison is with the combination of the Kolmogorov-Arnold model and the classical Gauss-Newton (GN) training method implemented in MATLAB by the authors of the present paper. The code is available via GitHub\footnote{https://github.com/mpoluektov/kan-polar}. In the implementation, the Tikhonov regularisation is used, $\lambda = 0.01$, and $50$ passes through the training dataset are made. The parameter updates are not scaled.

The second comparison is with the combination of the Kolmogorov-Arnold model and the Newton-Kaczmarz (NK) training method implemented in C\# by the authors of the present paper. The code is available via GitHub\footnote{https://github.com/andrewpolar/iai0}. In the implementation, $400$ passes through the training dataset are made. Furthermore, $\mu=0.1$ is used; however, the calculation of $\zeta$ is skipped for efficiency, taking $\zeta = 1$.

The final comparison is with a neural network model. The MATLAB's built-in neural network implementation, function \verb|fitrnet|, is chosen. For the benchmark, \verb|tanh| activation functions are chosen. Three layers, each of size $10$, are used (the default network configuration is one layer of size $10$), and $700$ passes through the training dataset are made. The philosophy behind constructing the architecture of the neural network is to use the most simple model that gives the desired accuracy, starting from the MATLAB's default. The model was obtained by increasing the number of layers to two and then to three (keeping the default size of $10$ for each layer), trying all MATLAB's activation functions meanwhile, and choosing the best. 

The training options for the methods implemented by the authors, as well as the number of passes for the neural network, have been attempted to be chosen such that the accuracy of the constructed models is close (apart from the GN method, which is difficult to fine-tune). However, the execution times are drastically different, with the C\# implementation of the NK method by the authors performing comparable to the MATLAB's built-in neural networks and around $14$ times faster than the Python implementation of \cite{Liu2024}. In this particular benchmark though, the GN method performs the best --- around $42$ times faster than the implementation of \cite{Liu2024}, while simultaneously being almost $3$ times more accurate. 

\subsubsection{Benchmark 2}

The benchmark above may create an illusion that the GN method in combination with the Kolmogorov-Arnold model with cubic-spline basis functions is the best possible option. This may not be true in general. For example, in section \ref{sec:examp1}, the GN method performed worse than the NK method.

To show that the NK method, as well as other basis functions, has a great potential, a more complex example is considered. It consists in fitting the function of $m=5$ inputs given by equation \eqref{eq:formula2}. The datasets consisting of $N = 10^4$ records for the training and of $N_\mathrm{val} = 1000$ records for the validation are used. The summary of the results is provided in table \ref{tab:times2}. The accuracy given in table \ref{tab:times2} is the normalised RMSE, i.e. $E_\mathrm{RMSE}$ that is given by equation \eqref{RMSE}; the values are given in $\%$.

\begin{table}
  \begin{center}
    \begin{tabular}{|l|ll|}
    \hline
    & accuracy ($\%$) & execution time (s) \\
    \hline
    present paper, C++ code, pw-l basis, NK method & $1.03 \pm 0.10$ & $0.055$ \\
    present paper, MATLAB code, spline basis, GN method & $1.01 \pm 0.23$ & $1.502$ \\
    MATLAB's neural network, \verb|relu| & $1.08 \pm 0.12$ & $0.265$ \\
    MATLAB's neural network, \verb|tanh| & $1.02 \pm 0.09$ & $0.408$ \\
    \hline
    \end{tabular}
  \end{center}
  \caption{The average normalised RMSEs and the average execution times for different methods and implementations --- benchmark 2. For the accuracy, the mean and the standard deviation across $10$ runs, are given.}
  \label{tab:times2}
\end{table}

In the C++ implementation by the authors, which is available via GitHub\footnote{https://github.com/andrewpolar/kanc2}, the `full' Kolmogorov-Arnold model ($11$ addends) with $n=6$ inner and $s=12$ outer piecewise-linear basis functions is constructed using the NK method; $36$ passes through the training dataset are made. In the implementation, $\mu=0.02$ is used for updating $H_{kjp}$ and $\mu=0.01$ is used for updating $G_{kl}$; the calculation of $\zeta$ is skipped for efficiency, taking $\zeta = 1$.

In the MATLAB implementation by the authors, the Kolmogorov-Arnold model of the same size as above (the number of addends and basis functions) is constructed, but with cubic-spline basis functions and using the GN method. The Tikhonov regularisation is used, $\lambda = 1$, and $13$ passes through the training dataset are made.

For the MATLAB's built-in neural network implementation, \verb|relu| and \verb|tanh| activation functions are used; three layers, each of size $10$, are used; $150$ and $100$ passes through the training dataset are made for the activation functions above, respectively. The same philosophy behind constructing the architecture of the neural network as in benchmark 1 is adopted.

Again, the training options for the methods are chosen such that the accuracy of the constructed models is close. In this benchmark, the GN method in combination with cubic-spline basis functions performs worse than the NK method in combination with piecewise-linear basis functions. The reason for the difference to benchmark 1 is the size of the problem --- the number of parameters\footnote{For the `full' model, the number of parameters is $\left(2m+1\right)\left(mn+s\right)$.} to be estimated ($462$ vs. $90$) and the number of training records ($10^4$ vs. $10^3$). Therefore, it may be expected that the GN method in combination with cubic-spline basis functions will not be the optimal choice for large models and/or datasets. In this example, the C++ implementation of the NK method by the authors performs $4$--$8$ times faster than the MATLAB's built-in neural networks.

\subsection{Additional examples: real-world datasets}

Although the scope of the paper is limited to the model's training method, there is a considerable interest of the ML community in performance of the Kolmogorov-Arnold model on real-world datasets and how the model compares to other models. Comprehensive investigation of this research question requires multiple separate studies. The aim of this subsection is to provide some contribution to this question. 

Results of five applications of the Kolmogorov-Arnold model are provided below, as well as its comparison to other models. In this subsection, only partial details are provided to avoid diluting the main text of the present paper, and all simulation details can be obtained from the corresponding open-source codes.

The first dataset contains results of aerodynamic and acoustic tests of 2D and 3D airfoil blade sections\footnote{https://archive.ics.uci.edu/dataset/291/airfoil+self+noise}. It has $1503$ records and $5$ inputs. The accuracy metric is the Pearson correlation coefficient, which was $0.95$ for the Kolmogorov-Arnold model\footnote{http://openkan.org/airfoil.html}. Third-party software Neural Designer achieved the same accuracy with the neural network model\footnote{https://www.neuraldesigner.com/learning/examples/airfoil-self-noise-prediction}.

The second dataset contains clients' decisions regarding a subscription to a bank service\footnote{https://www.kaggle.com/datasets/gauravtopre/bank-customer-churn-dataset}. It has $10000$ records and $11$ inputs. The accuracy metric is the fraction of correct predictions, which was $0.85$ for the Kolmogorov-Arnold model\footnote{http://openkan.org/bankchurn.html}. Third-party software Neural Designer achieved $0.79$ with the neural network model\footnote{https://www.neuraldesigner.com/learning/examples/bank-churn}.

The third dataset contains characteristics of medical images for cancer diagnostics\footnote{https://archive.ics.uci.edu/dataset/17/breast+cancer+wisconsin+diagnostic}. It has $569$ records and $30$ inputs. The accuracy metric is the fraction of correct predictions, which was $0.97$ for the Kolmogorov-Arnold model\footnote{http://openkan.org/wdbc.html}. The same accuracy was achieved using neural networks in \cite{Poeta2024}.

The fourth dataset contains characteristics of emails classified as spam and non-spam\footnote{https://archive.ics.uci.edu/dataset/94/spambase}. The accuracy metric is the fraction of correct predictions, which was $0.93$ for the Kolmogorov-Arnold model\footnote{http://openkan.org/spambase.html}. The accuracy of $0.94$ was achieved using neural networks in \cite{Poeta2024}.

The fifth dataset contains information on annual income of individuals\footnote{https://archive.ics.uci.edu/dataset/2/adult}. The accuracy metric is the fraction of correct predictions, which was $0.86$ for the Kolmogorov-Arnold model\footnote{http://openkan.org/adult.html}. The same accuracy was achieved using neural networks in \cite{Poeta2024}.

\subsection{Large-scale example}

The purpose of this subsection is to compare further the Kolmogorov-Arnold model with neural networks, but now using a large-scale example for a high-performance laptop --- an example that requires hours of training time.

The chosen example is constructing a regression model for a synthetic dataset --- the input is a $5$-by-$5$ matrix and the output is the determinant of the matrix:
\begin{equation*}
  y_i = \begin{vmatrix} x_{i,1} & \ldots & x_{i,5} \\
    \vdots & \ddots & \vdots \\
    x_{i,21} & \ldots & x_{i,25} \end{vmatrix} .
\end{equation*}
The datasets are generated by taking uniformly-distributed random inputs $x_{i,j} \sim \operatorname{unif}\left(0,1\right)$ and calculating outputs $y_i$.

It transpires that training a regression model for such problem is a relatively difficult task due to a very high non-linearity and requires large datasets (even though problems with matrices up to $4$-by-$4$ are relative easy, with small datasets being sufficient). Therefore, the training and the validation datasets' sizes are taken to be $10^7$ and $2 \cdot 10^6$ records, respectively. The performance of the trained models is assessed by calculating both the normalised RMSE and the Pearson correlation coefficient between the original output and the modelled output on the validation datasets. 

The state-of-the-art implementation of neural networks is chosen for the comparison --- the MATLAB's built-in neural network implementation. For the benchmark, \verb|tanh| activation functions are chosen and the networks consisting of the fully-connected layers are trained via function \verb|fitrnet|. For simplicity, the number of neurons is kept the same for all layers. The number of layers and the number of neurons per layer are varied, as well as the total number of iterations.

The C++ implementation by the authors of the Kolmogorov-Arnold model in combination with the Newton-Kaczmarz (NK) training method is used for the comparison. The piecewise-linear basis functions are taken. The number of addends, as well as the number of inner and outer basis functions, is varied. 

\begin{figure}
  \begin{center}
    \includegraphics[scale=0.9]{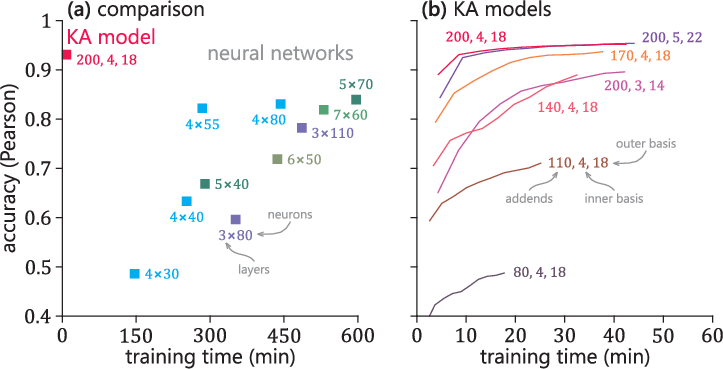}
  \end{center}
  \caption{The comparison of the performance of the Kolmogorov-Arnold (KA) model and the MATLAB's built-in neural networks (a) in the problem of building a regression model of a determinant of a $5$-by-$5$ matrix. Evolution of the accuracy of the KA model (measured on the validation dataset) after each pass through the training dataset (b).}
  \label{fig:det5}
\end{figure}

The results are summarised in figure \ref{fig:det5} and table \ref{tab:compNew}. In this example, the Kolmogorov-Arnold model significantly outperforms the neural networks both in terms of accuracy (measured via the Pearson correlation coefficient) and training time. The `best' configuration of the neural network ($4$ layers with $55$ neurons per layer) gave the correlation of $0.822$ after $4$ hours $44$ minutes of training, requiring $1000$ iterations, while all tested configurations are shown in figure \ref{fig:det5}a. The `best' configuration of the Kolmogorov-Arnold model ($200$ addends, $4$ inner basis functions, $18$ outer basis functions) gave the correlation of $0.890$ only after $1$ pass through the dataset, running for $4$ minutes $15$ seconds. Further training (passes through the training dataset) improves the accuracy of the Kolmogorov-Arnold model, as shown in figure \ref{fig:det5}b, where the curves are the accuracy calculations on the validation dataset after each pass through the training dataset.

\begin{table}
  \begin{center}
    \begin{tabular}{|l|l|lllll|}
    \hline
    & configuration & Pearson & RMSE & runtime (min) & CPU & memory ($\text{GB}$) \\
    \hline
    \multirow{7}{*}{KA} & $200,5,22$; $1$ pass & $0.843$ & $0.018$ & $4.66$ & $9$--$11\%$ & $2.6$ \\
    & $200,5,22$; $2$ passes & $0.925$ & $0.013$ & $9.28$ & $9$--$11\%$ & $2.6$ \\
    & $200,5,22$; $10$ passes & $0.954$ & $0.010$ & $44.17$ & $9$--$11\%$ & $2.6$ \\
    & $200,4,18$; $1$ pass & $0.890$ & $0.016$ & $4.26$ & & \\
    & $200,4,18$; $2$ passes & $0.931$ & $0.013$ & $8.51$ & & \\
    & $170,4,18$; $10$ passes & $0.937$ & $0.013$ & $37.78$ & & \\
    & $140,4,18$; $10$ passes & $0.890$ & $0.015$ & $32.58$ & & \\
    \hline
    \multirow{7}{*}{NN} & $4 \!\times\! 55$; $1000$ iter. & $0.822$ & $0.019$ & $284.15$ & & \\
    &  $4 \!\times\! 80$; $1000$ iter. & $0.831$ & $0.018$ & $443.08$ & $33$--$35\%$ & $45.7$ \\
    &  $5 \!\times\! 70$; $1000$ iter. & $0.840$ & $0.018$ & $596.50$ & $33$--$35\%$ & $47.7$ \\
    &  $7 \!\times\! 60$; $1000$ iter. & $0.819$ & $0.019$ & $530.63$ & $33$--$35\%$ & $50.4$ \\
    &  $4 \!\times\! 80$; $100$ iter. & $0.002$ & $0.033$ & $46.81$ & & \\
    &  $5 \!\times\! 70$; $100$ iter. & $0.007$ & $0.033$ & $55.16$ & & \\
    &  $7 \!\times\! 60$; $100$ iter. & $0.007$ & $0.033$ & $58.02$ & & \\
    \hline
    \end{tabular}
  \end{center}
  \caption{The summary of some results of training of the Kolmogorov-Arnold (KA) model and the MATLAB's built-in neural networks (NN). Various configurations are used: for KA, the numbers indicate the number of addends, the number of inner basis functions, the number of outer basis functions, the number of completed passes through the training dataset; for NN, the numbers indicate the number of layers, the number of neurons per each layer, the number of completed iterations utilising the training dataset.}
  \label{tab:compNew}
\end{table}

As shown in table \ref{tab:compNew}, the neural networks indeed require large number of iterations to achieve acceptable accuracy. After $100$ iterations, taking between $45$ minutes and $1$ hour of training, the correlation is still around $0$, making such models unusable. Meanwhile, $45$ minutes of training was sufficient for the Kolmogorov-Arnold model to achieve the ultimate correlation of $0.954$, as shown in figure \ref{fig:det5}b.

Both implementations were executed on a Dell XPS laptop with Intel Core i7-12700H processor and $64\,\text{GB}$ of RAM. The latter is essential, as running the neural-network training was taking up to $50\,\text{GB}$ of RAM. Furthermore, the MATLAB's implementation of neural networks uses several cores with the default value of \verb|maxNumCompThreads| being $6$ on the particular system. The CPU has $6$ \emph{performance} and $8$ \emph{efficiency} cores, and the CPU usage for neural network training was reported to be in the range of $33$--$35\%$, indicating that probably only $4$ cores were used. The non-parallel C++ implementation the Kolmogorov-Arnold model was used, hence, it utilised strictly $1$ core. Furthermore, its memory footprint was more than an order of magnitude lower --- below $3\,\text{GB}$. Such drastic difference in memory footprints is due to the NK method not requiring assembling of large matrices; the only matrices that are held in memory are the datasets (training and validation) and the parameters.

The codes are available online and the results can be easily reproduced. The script constructing the neural networks is elementary and can be found on the authors' website\footnote{http://openkan.org/tetrahedron.html}. The C++ code is available via GitHub\footnote{https://github.com/andrewpolar/det5}.

\section{Conclusions}
\label{sec:conclusion}

The Kolmogorov-Arnold model is a convenient tool for tasks that require mapping input vectors of a black box system into output scalars. In the present paper, it has been proposed to decompose the underlying functions constituting the model into the basis functions and the model parameters, and to determine these parameters using the Newton-Kaczmarz (NK) method, which is locally convergent.

The approach has been tested using the numerical examples. First, it has been compared to the straightforward application of the Gauss-Newton method to estimate the parameters of the ridge function model. A challenging example has been selected, where the non-linearity of the outer function creates convergence difficulties when the initial guess for the parameters is far from their exact values. It has been shown that the NK method is much less sensitive to the selection of the initial guess; therefore, it is more robust from a practical point of view. The NK method has also been compared to the `Adam' stochastic gradient descent method (SGD), and it has been shown that although it is possible to get the same performance out of both methods, it is subjectively easier to fine-tune the NK method to achieve the best performance\footnote{In particular, the `Adam' SGD method has three fine-tuning parameters, while the NK method has only one.}.

The second example has demonstrated that a challenging non-linear function can be approximated accurately by the Kolmogorov-Arnold model with the piecewise-linear basis functions, and the model parameters can be obtained efficiently using the proposed approach. 

Furthermore, it has been shown that the Kolmogorov-Arnold model can be used for data-driven solution of partial differential equations (PDEs). There is crucial difference between the proposed approach and employing neural networks for such task --- the underlying functions of the Kolmogorov-Arnold model are taken as products of the model parameters and the basis functions. The choice of the basis functions with compact support, similarly to e.g. finite-element methods, can be viewed as discretisation of the spatial and the temporal coordinates --- the feature that physics-informed neural networks (PINNs) lack. Thus, the resulting approach can be viewed as a hybrid between the data-driven and the mesh-based solvers.

Model training using the NK method is ideally-suited for parallel implementation. The present paper has shown how the parameter update and the basis function evaluation can be decoupled, leading to an option for multi-threaded implementation. Preliminary tests show that such modification does not result in a noticeable loss of accuracy. A possibility of such parallel implementation is a significant advantage of the proposed approach compared to other parameter estimation methods for the Kolmogorov-Arnold model.

Finally, since multiple implementations by the authors have been discussed, it is useful to mention that the authors recommend their latest MATLAB implementation\footnote{https://github.com/mpoluektov/kan-polar} and their latest C++ implementation\footnote{https://github.com/andrewpolar/kanpro2}, depending on the programming language preference of the user, as good starting points.

\section*{Competing interests}

The authors have no conflicts of interest to declare that are relevant to the content of this article.

\appendix
\titleformat{\section}[hang]{\Large\bfseries\raggedright\sffamily}{Appendix \thesection}{1em}{}

\section{The Newton-Kaczmarz method}
\label{sec:NKmethod}

For system of non-linear equations $\vn{L}\left(\vn{Z}\right) = 0$, where $\vn{L}: D \to \mathbb{R}^N$ is continuously differentiable in open set $D \subset \mathbb{R}^C$, the Newton-Kaczmarz method \cite{Meyn1983,MARTINEZ1986,MARTINEZ1986a} consists in iterations
\begin{equation}
  \vn{Z}^{q+1} = \vn{Z}^q - \mu\frac{L_i\left(\vn{Z}^q\right)}{\left\| \nabla L_i\left(\vn{Z}^q\right) \right\|^2} {\nabla L_i\left(\vn{Z}^q\right)} , \quad\quad
  \nabla L_i\left(\vn{Z}^q\right) = \left. \begin{bmatrix} \frac{\partial L_i}{\partial Z_1} & \!\!\ldots\!\! & \frac{\partial L_i}{\partial Z_C} \end{bmatrix}^\mathrm{T} \right|_{\vn{Z} = \vn{Z}^q} ,
  \label{eq:NewtKacz}
\end{equation}
where function $L_i$ denotes the $i$-th component of vector function $\vn{L}$, scalar $Z_j$ denotes the $j$-th component of vector $\vn{Z}$, parameter $\mu \in \left(0,1\right]$ is introduced for numerical damping, and vector $\vn{Z}^q$ is the approximation of the solution at iteration $q$. Index $i$ changes each iteration. The method requires an initial guess and continues until some convergence criteria are reached. 

This method can also be interpreted as an iterative scheme, where at each step, a single equation is linearised and one step of the Kaczmarz method \cite{Kaczmarz1937,Tewarson1969} is performed. Indeed, to write an iterative scheme, an expression for $\vn{Z}^{q+1}$ must be found, assuming that $\vn{Z}^q$ is known. Considering single equation $L_i(\vn{Z}^{q+1}) = 0$, the linearisation is performed:
\begin{equation}
  L_i(\vn{Z}^{q+1}) = L_i\left(\vn{Z}^q\right) + \nabla L_i\left(\vn{Z}^q\right)^\mathrm{T} \Delta \vn{Z} + O\left(\left\|\Delta \vn{Z}\right\|^2\right) = 0 , \quad\quad 
  \Delta \vn{Z} = \vn{Z}^{q+1} - \vn{Z}^q .
\end{equation}
When the second-order term is neglected, linear equation 
\begin{equation}
  \vn{M}_i \Delta \vn{Z} = Y_i 
  \label{eq:MdltZ}
\end{equation}
is obtained, where $\vn{M}_i = \nabla L_i\left(\vn{Z}^q\right)^\mathrm{T}$ and $Y_i = -L_i\left(\vn{Z}^q\right)$. Equation \eqref{eq:MdltZ} with respect to $\Delta \vn{Z}$ describes a hyperplane in $\mathbb{R}^C$. Using the key idea of the Kaczmarz method (i.e. the projection descend), the origin is projected onto this hyperplane to obtain the solution:
\begin{equation}
  \Delta \vn{Z} = \frac{Y_i}{\left\| \vn{M}_i \right\|^2} \vn{M}_i^\mathrm{T} .
\end{equation}
Adding numerical damping parameter $\mu$ to limit the magnitude of change of $\vn{Z}$ over one step, leads to equation \eqref{eq:NewtKacz}.

\subsection{On local convergence of the Newton-Kaczmarz method}
\label{sec:errAn}

As the Newton-Kaczmarz method is a version of the Newton's method, it is not expected to be globally convergent (apart from isolated cases, e.g. when $\vn{L}\left(\vn{Z}\right)$ are linear equations), but it is locally convergent, given a `good' initial guess. It is easy to see this by performing a simple error analysis. 

System of equations $\vn{L}\left(\vn{Z}\right) = 0$ is assumed to be consistent, i.e. there exists at least one $\bar{\vn{Z}}$ that satisfies all equations. It is now assumed that $\vn{L}: D \to \mathbb{R}^N$ is twice continuously differentiable in open set $D \subset \mathbb{R}^C$. Initial guess $\vn{Z}^0$ is assumed to be sufficiently close to $\bar{\vn{Z}}$. Using the Taylor's expansion results in
\begin{equation}
  L_i(\bar{\vn{Z}}) = -Y_i + \vn{M}_i \vn{E}^q + R_i = 0 , \quad\quad\quad
  \vn{E}^q = \bar{\vn{Z}} - \vn{Z}^q , \quad\quad\quad
  R_i = O\left(\left\|\vn{E}^q\right\|^2\right) ,
  \label{eq:Ldecom}
\end{equation}
where the notation introduced in equation \eqref{eq:MdltZ} is utilised. Using iterative formula \eqref{eq:NewtKacz} and subtracting $\bar{Z}$ from both sides results in
\begin{equation}
  -\vn{E}^{q+1} = -\vn{E}^q + \mu\frac{Y_i}{\left\| \vn{M}_i \right\|^2} \vn{M}_i^\mathrm{T} .
  \label{eq:errit}
\end{equation}
It is useful to remind that for any vector $\vn{A}$, the norm is defined as $\left\| \vn{A} \right\|^2 = \vn{A}^\mathrm{T} \vn{A}$, if $\vn{A}$ is a column, and $\left\| \vn{A} \right\|^2 = \vn{A} \vn{A}^\mathrm{T}$, if $\vn{A}$ is a row. Equation \eqref{eq:errit} then gives
\begin{equation}
\begin{split}
  \left\| \vn{E}^{q+1} \right\|^2 &= \left\| \vn{E}^q \right\|^2 + \mu^2\frac{Y_i^2}{\left\| \vn{M}_i \right\|^2} - 2 \mu\frac{Y_i}{\left\| \vn{M}_i \right\|^2} \vn{M}_i \vn{E}^q = \\
  &= \left\| \vn{E}^q \right\|^2 - \frac{\mu}{\left\| \vn{M}_i \right\|^2}\left(\left(2-\mu\right) Y_i^2 - 2 Y_i R_i \right) ,
\end{split}
  \label{eq:errdec}
\end{equation}
where $\vn{M}_i \vn{E}^q$ is substituted from \eqref{eq:Ldecom}. To have local convergence, it is sufficient that $\left\| \vn{E}^{q+1} \right\|^2 < \left\| \vn{E}^q \right\|^2$, which means that the error decreases at each iteration. Hence, using equation \eqref{eq:errdec}, it remains to prove that
\begin{equation*}
  \left(2-\mu\right) Y_i^2 - 2 Y_i R_i > 0 .
\end{equation*}
If $0 < \mu \leq 1$, then 
\begin{equation}
  \left(2-\mu\right) Y_i^2 - 2 Y_i R_i \geq Y_i^2 - 2 Y_i R_i = 
  \left( \vn{M}_i \vn{E}^q \right)^2 - R_i^2 ,
  \label{eq:ineqConv}
\end{equation}
where $Y_i$ is substituted from \eqref{eq:Ldecom}. Assumption that $\nabla L_i$ is finite in the neighbourhood of solution $\bar{\vn{Z}}$ leads to 
\begin{equation*}
  \left( \vn{M}_i \vn{E}^q \right)^2 = O\left(\left\|\vn{E}^q\right\|^2\right) .
\end{equation*}
The first and the second terms of the right-hand side of equation \eqref{eq:ineqConv} become second-order and fourth-order, respectively, with respect to $\left\|\vn{E}^q\right\|$. The second-order term is positive, while the fourth-order term is negative. For small $\left\|\vn{E}^q\right\|$, the second-order term is dominant; hence, the right-hand side of equation \eqref{eq:ineqConv} is positive when $\vn{Z}^q$ is close to $\bar{\vn{Z}}$. The latter means that the iterative scheme is locally convergent.

It is useful to note that in equation \eqref{eq:errdec}, term $Y_i^2$ is the squared residual of $i$-th equation with the solution taken at $q$-th iteration. It is the leading-order term and has multiplier $\left(2\mu-\mu^2\right)$. Thus, neglecting the higher-order terms containing Taylor remainder $R_i$, the fastest error decay is achieved when this multiplier takes the maximum value, i.e. when $\mu = 1$. The purpose of introducing $\mu$ within the iterative scheme is to make the method robust for noisy data, i.e. when system of equations $\vn{L}\left(\vn{Z}\right) = 0$ is not consistent, in which case, the decrease of $\mu$ acts as noise filtering.

\subsection{On relation to the stochastic gradient descent method}
\label{sec:SGD}

It can be useful to note that iterative formula \eqref{eq:NewtKacz} has the same structure as the well-known stochastic gradient descent (SGD) method. The SGD method is written for minimisation of objective function 
\begin{equation*}
  Q \left( \vn{Z} \right) = \frac{1}{N} \sum_{i=1}^N Q_i \left( \vn{Z} \right) ,
\end{equation*}
where $\vn{Z}$ is the vector of parameters to be found by the minimisation of $Q \left( \vn{Z} \right)$. Addends $Q_i \left( \vn{Z} \right)$ are typically associated with the $i$-th record of the dataset. The SGD method is then written as the following iterative formula:
\begin{equation*}
  \vn{Z}^{q+1} = \vn{Z}^q - \tilde{\mu} \, \nabla Q_i\left( \vn{Z}^q \right) ,
\end{equation*}
where $\tilde{\mu}$ is referred to as the ``learning rate''. Index $i$ changes each iteration. Iterative formula \eqref{eq:NewtKacz} is obtained by defining 
\begin{equation*}
  Q_i\left( \vn{Z} \right) = \frac{1}{2} \left(L_i\left( \vn{Z} \right)\right)^2 , \quad\quad\quad 
  \tilde{\mu} = \mu \frac{1}{\left\| \nabla L_i\left(\vn{Z}^q\right) \right\|^2} .
\end{equation*}
Hence, the iterations of the Newton-Kaczmarz method can be viewed as the iterations of the SGD method with an adaptive learning rate. It should be emphasised that the convergence analysis of the SGD method with convex and non-convex objective functions is well studied, e.g. \cite{Fehrman2020}.

\section{Derivatives of the Kolmogorov-Arnold model}
\label{sec:derv}

The implementation of the Newton-Kaczmarz method to find the parameters of the Kolmogorov-Arnold model requires expressions for the derivatives of the output by the parameters. Furthermore, solving PDEs using the proposed approach additionally requires expressions for the mixed derivatives of the output by the parameters and the inputs. These derivatives are summarised in this section. The differentiation with respect to the inputs is limited to the second order.  

It is useful to start by writing the Kolmogorov-Arnold model with the substituted basis functions --- substituting equations \eqref{eq:KAfunc1} and \eqref{eq:KAfunc2} into \eqref{eq:KAU}:
\begin{equation}
  \tilde{y} = \sum_{k,l} G_{kl} \psi_l \left( \sum_{j,p} H_{kjp} \phi_p \left(x_j\right) \right) ,
  \label{eq:KASubs}
\end{equation}
where superscript $i$ is omitted for brevity, as well as the ranges of the indices used in the summations. The following short notation is introduced for the basis functions evaluated at the inputs and at the intermediate variables:
\begin{equation*}
  \varLambda_{pj} = \phi_p \left(x_j\right) ,
  \quad\quad\quad
  \varPsi_{lk} = \psi_l \left( \sum_{j,p} H_{kjp} \varLambda_{pj} \right) .
\end{equation*}
The derivatives of the basis functions evaluated at the same points are denoted with the prime symbol, for example,
\begin{equation*}
  \varLambda_{pj}' = \left.\frac{\dif \phi_p\left(x\right)}{\dif x}\right|_{x=x_j} . 
\end{equation*}
The derivatives of $\tilde{y}$ by the parameters and by the inputs are the following:
\begin{align*}
  &\frac{\partial \tilde{y}}{\partial G_{kl}} = \varPsi_{lk} , \\
  &\frac{\partial \tilde{y}}{\partial H_{kjp}} = \sum_{l} G_{kl} \varPsi_{lk}' \varLambda_{pj} , \\
  &\frac{\partial \tilde{y}}{\partial x_a} = \sum_{k,l,p} G_{kl} H_{kap} \varPsi_{lk}' \varLambda_{pa}' .
\end{align*}
The second-order derivatives are the following:
\begin{align*}
  &\frac{\partial^2 \tilde{y}}{\partial x_a \partial G_{kl}} = \sum_{p} H_{kap} \varPsi_{lk}' \varLambda_{pa}' , \\
  &\frac{\partial^2 \tilde{y}}{\partial x_a \partial H_{kjp}} = \sum_{l} G_{kl} \delta_{aj} \varPsi_{lk}' \varLambda_{pa}' + 
  \sum_{l,u} G_{kl} H_{kau} \varPsi_{lk}'' \varLambda_{ua}' \varLambda_{pj} , \\
  &\frac{\partial^2 \tilde{y}}{\partial x_a \partial x_b} = \sum_{k,l,p,u} G_{kl} H_{kap} H_{kbu} \varPsi_{lk}'' \varLambda_{pa}' \varLambda_{ub}' +
  \sum_{k,l,p} G_{kl} H_{kap} \delta_{ab} \varPsi_{lk}' \varLambda_{pa}'' ,
\end{align*}
where $\delta_{ab}$ is the Kronecker delta. The third-order derivatives are the following:
\begin{align*}
  &\frac{\partial^3 \tilde{y}}{\partial x_a \partial x_b \partial G_{kl}} = \sum_{p,u} H_{kap} H_{kbu} \varPsi_{lk}'' \varLambda_{pa}' \varLambda_{ub}' + \sum_{p} H_{kap} \delta_{ab} \varPsi_{lk}' \varLambda_{pa}'' , \\
  &\frac{\partial^3 \tilde{y}}{\partial x_a \partial x_b \partial H_{kjp}} = \sum_{l,u} G_{kl} H_{kbu} \delta_{aj} \varPsi_{lk}'' \varLambda_{pa}' \varLambda_{ub}' + \sum_{l,u} G_{kl} H_{kau} \delta_{bj} \varPsi_{lk}'' \varLambda_{pb}' \varLambda_{ua}' + \\ &\hspace{0.5cm}+
  \sum_{l,u,v} G_{kl} H_{kbu} H_{kav} \varPsi_{lk}''' \varLambda_{va}' \varLambda_{ub}' \varLambda_{pj} + \sum_{l} G_{kl} \delta_{aj} \delta_{ab} \varPsi_{lk}' \varLambda_{pa}'' +
  \sum_{l,u} G_{kl} H_{kau} \delta_{ab} \varPsi_{lk}'' \varLambda_{ua}'' \varLambda_{pj} .
\end{align*}

\section{Parameter estimation of the Kolmogorov-Arnold model using the Gauss-Newton method}
\label{sec:app}

The Kolmogorov-Arnold model with the substituted basis functions is given by equation \eqref{eq:KASubs}. The sum of squares of the residuals is introduced as
\begin{equation}
  S = \frac{1}{2} \sum_{i=1}^{N} \left( \tilde{y}_i - y_i \right)^2 .
\end{equation}
The Gauss-Newton method consists in minimisation of $S$ with respect to the model parameters using the Newton-Raphson method and neglecting the second-order derivative terms in the Jacobian. The minimisation of $S$ gives
\begin{equation}
  \frac{\partial S}{\partial G_{kl}} = \sum_i r_i \frac{\partial \tilde{y}}{\partial G_{kl}} = 0 , \quad\quad
  \frac{\partial S}{\partial H_{kjp}} = \sum_i r_i \frac{\partial \tilde{y}}{\partial H_{kjp}} = 0 , 
  \quad\quad r_i = \tilde{y}_i - y_i ,
  \label{eq:ridgeSetEqs}
\end{equation}
where, as previously, it is implied that values of the derivatives of $\tilde{y}$ are calculated at $\vn{x}_i$. The evaluation of the second derivatives of $S$ leads to 
\begin{align}
  &\frac{\partial^2 S}{\partial G_{kl} \partial G_{ab}} = \sum_i \frac{\partial \tilde{y}}{\partial G_{kl}} \frac{\partial \tilde{y}}{\partial G_{ab}} , \quad\quad\quad
  \frac{\partial^2 S}{\partial G_{kl} \partial H_{ajp}} = \sum_i \left( \frac{\partial \tilde{y}}{\partial G_{kl}} \frac{\partial \tilde{y}}{\partial H_{ajp}} + r_i \frac{\partial^2 \tilde{y}}{\partial G_{kl} \partial H_{ajp}} \right) , \nonumber \\
  &\frac{\partial^2 S}{\partial H_{kjp} \partial H_{abu}} = \sum_i \left( \frac{\partial \tilde{y}}{\partial H_{kjp}} \frac{\partial \tilde{y}}{\partial H_{abu}} + r_i \frac{\partial^2 \tilde{y}}{\partial H_{kjp} \partial H_{abu}} \right) , 
  \label{eq:dir2S}
\end{align}
Using the obtained relations, equations \eqref{eq:ridgeSetEqs} are solved using the Newton-Raphson method, consisting in the following iterative scheme:
\begin{equation}
  \vn{Z}^{q+1} = \vn{Z}^q - \mu \left. \left(\vn{M} + \lambda\vn{I}\right)^{-1} \vn{L} \right|_{\vn{Z}=\vn{Z}^q} ,
  \label{eq:GNupd}
\end{equation}
where $\vn{Z}$ is the vector containing the model parameters, with superscript indicating the iteration number, $\vn{M}$ and $\vn{L}$ are the matrix and the vector containing the second and the first derivatives of $S$ with respect to the model parameters, respectively, $\lambda$ and $\vn{I}$ are the Tikhonov regularisation parameter and the identity matrix, respectively, and $\mu$ is the numerical damping parameter added for the robustness of the scheme. Vectors $\vn{Z}$ and $\vn{L}$ can be formally written as
\begin{align}
  &\vn{Z} = \begin{bmatrix}
    G_{11} & \!\!\ldots\!\! & G_{ds} & H_{111} & \!\!\ldots\!\! & H_{dmn}
  \end{bmatrix}\transp , \\
  &\vn{L} = \begin{bmatrix}
    \frac{\partial S}{\partial G_{11}} & \!\!\ldots\!\! & \frac{\partial S}{\partial G_{ds}} & \frac{\partial S}{\partial H_{111}} & \!\!\ldots\!\! & \frac{\partial S}{\partial H_{dmn}}
  \end{bmatrix}\transp ,
\end{align}
with matrix $\vn{M}$ assembled accordingly from derivatives \eqref{eq:dir2S}, neglecting the second derivatives of $\tilde{y}$, and is omitted here for brevity. The iterations are performed until either the maximum number of iterations is reached or $\| \vn{Z}^{q+1} - \vn{Z}^q \| < \delta$ is fulfilled.

\section{Summary of the `Adam' SGD method}
\label{sec:adam}

The `Adam' SGD method has been proposed in \cite{Kingma2015}. Using the notation of appendix \ref{sec:SGD}, the method constitutes the following iterative sequence:
\begin{align*}
  &\vn{G}^q = \nabla Q_i\left( \vn{Z}^q \right) , \\
  &\vn{A}^{q+1} = \beta_1 \vn{A}^q + \left(1-\beta_1\right) \vn{G}^q , \quad\quad\quad
  B^{q+1}_j = \beta_2 B^q_j + \left(1-\beta_2\right) \left(G^q_j\right)^2 , \quad\quad\quad j \in \left\lbrace 1, \ldots, C \right\rbrace , \\
  &\hat{\vn{A}}^{q+1} = \frac{\vn{A}^{q+1}}{1 - \left(\beta_1\right)^q} , \quad\quad\quad\quad
  \hat{\vn{B}}^{q+1} = \frac{\vn{B}^{q+1}}{1 - \left(\beta_2\right)^q} , \quad\quad\quad\quad
  Z^{q+1}_j = Z^q_j - \mu \frac{\hat{A}^{q+1}_j}{\sqrt{\hat{B}^{q+1}_j} + \varepsilon} ,
\end{align*}
where vector $\vn{Z}^q \in \mathbb{R}^C$ is the approximation of the sought solution at iteration $q$, vector $\vn{G}^q \in \mathbb{R}^C$ is the gradient of objective function $Q_i$ (index $i$ changes each iteration), vectors $\vn{A}^q$ and $\vn{B}^q$ are referred to as the biased first and second moment estimates, respectively, vectors $\hat{\vn{A}}^q$ and $\hat{\vn{B}}^q$ are the corresponding bias-corrected moments, parameters $\beta_1,\beta_2 \in \left[0,1\right)$ are referred to as the `forgetting' factors, small scalar parameter $\varepsilon$ is introduced to avoid division by zero. Scalar $B^q_j$ denotes the $j$-th component of vector $\vn{B}^q$ (the same for $\vn{A}^q$ and $\vn{G}^q$). The initial values for $\vn{A}^q$ and $\vn{B}^q$ are the all-zero vectors.

\bibliographystyle{unsrt}
\bibliography{refs}

\end{document}